\documentclass[12pt]{amsart}
\usepackage{latexsym, amssymb, amsmath, amscd}
 \usepackage{eucal}

    \newtheorem{rema}{Remark}[section]
    \newtheorem{propo}[rema]{Proposition}

   \newtheorem{theo}[rema]{Theorem}
   \newtheorem{def-theo}[rema]{Definition-Theorem}
 
   \newtheorem{defi}[rema]{Definition}
    \newtheorem{lemma}[rema]{Lemma}
    \newtheorem{corol}[rema]{Corollary}
     \newtheorem{exam}[rema]{Example}
  \newtheorem{rmk}[rema]{Remark}

	\newcommand{\nno}{\nonumber}

\renewcommand{\hom}{\mbox{\rm Hom}}
 \newcommand{\pf}{{\it Proof:}\hspace{2ex}}
 
 \newcommand{\epfv}{\hspace{1em}$\Box$\vspace{1em}}

\newcommand{\bC}{{\mathbb C}}
\newcommand{\bZ}{{\mathbb Z}}
\newcommand{\bR}{{\mathbb R}}

\newcommand{\bN}{{\mathbb N}}

\newcommand{\cA}{{\mathcal A}}

\newcommand{\Tr}{\mbox{\rm Tr\,}}

\newcommand{\rad}{{\frak r}\hspace{0.3mm}}
\newcommand{\arad}{ {\frak r\hspace{0.3mm}'} }
\newcommand{\nil}{\mbox{\rm nil\,}}

\newcommand{\vms}{$\vartheta$-Mathieu subspace }
\newcommand{\vmss}{$\vartheta$-Mathieu subspaces }
\newcommand{\gca}{\mathcal{G}(\cA)}
\newcommand{\ovca}{\mathcal{E}_\vartheta(\cA)}

\newcommand{\cB}{{\mathcal B}}

\newcommand{\cF}{{\mathcal F}}
\newcommand{\cE}{{\mathcal E}}

\newcommand{\Ker}{\mbox{\rm Ker\,}}

\newcommand{\vt}{\vartheta}

\title[Mathieu Subspaces of Associative Algebras]
{Mathieu Subspaces of Associative Algebras}

  \author{Wenhua Zhao}      
\address{Department of Mathematics, Illinois State University, Normal, IL 61761. Email: wzhao@ilstu.edu}

\begin{document}

\begin{abstract}
Motivated by {\it the Mathieu conjecture} \cite{Ma}, 
 {\it the image conjecture} \cite{IC} and 
the well-known {\it Jacobian conjecture} \cite{K} 
(see also \cite{BCW} and \cite{E}), the notion of 
{\it Mathieu subspaces} as a natural  generalization 
of the notion of ideals has been introduced recently 
in \cite{GIC} for associative algebras. In this paper, 
we first study algebraic elements in the radicals of 
Mathieu subspaces of associative algebras over fields 
and prove some properties and characterizations 
of Mathieu subspaces with algebraic radicals. 
We then give some characterizations or classifications for 
{\it strongly simple algebras} (the algebras 
with no non-trivial  Mathieu subspaces)  
over arbitrary commutative rings, 
and for {\it quasi-stable algebras}   
(the algebras all of whose subspaces 
that do not contain the identity element 
of the algebra are Mathieu subspaces)  
over arbitrary fields. 
Furthermore, co-dimension 
one Mathieu subspaces and the minimal non-trivial  
Mathieu subspaces of the matrix algebras over fields 
are also completely determined.  
\end{abstract}

   

%16N40 Nil and nilpotent radicals, sets, ideals, rings
%
%16D70 Structure and classification (except as in 16Gxx), direct sum decomposition, cancellation
%
%16D99 None of the above, but in this section

%16S34 Group rings [see also 20C05, 20C07], Laurent polynomial rings

%13N10: Rings of differential operators and their modules 
%14R15: Jacobian problem.
%33C45 Orthogonal polynomials and functions of hypergeometric type (Jacobi, Laguerre, Hermite, Askey scheme, etc.) [see also 42C05 for general orthogonal polynomials and functions]. 
%32C38 Sheaves of differential operators and their modules, $D$-modules. 
%32W99: Differential operators in several 
%    complex variables/None of the above, but in this section. 

\thanks{The author has been partially supported 
by NSA Grant H98230-10-1-0168}

 \bibliographystyle{alpha}
    \maketitle

 % \tableofcontents

\renewcommand{\theequation}{\thesection.\arabic{equation}}
\renewcommand{\therema}{\thesection.\arabic{rema}}
\setcounter{equation}{0}
\setcounter{rema}{0}
\setcounter{section}{0}

\section{\bf Introduction}

\subsection{Background and Motivation}
Let $R$ be an arbitrary commutative ring and $\cA$ 
an associative but not necessarily commutative 
algebra over $R$. Then we have the following notion 
introduced recently by the author in \cite{GIC}.
 
\begin{defi}\label{Def-MS}
Let $M$ be a  $R$-submodule or $R$-subspace of $\cA$. 
We say $M$ is a {\it left} $($resp.,  {\it right}$)$ 
{\it  Mathieu subspace} of $\cA$  
if the following property holds: 
let $a\in \cA$ such that $a^m\in M$ for all $m\ge 1$. 
Then for any $b\in \cA$, we have 
$ b a^m \in M$ $($resp.,  $a^m b \in M$$)$ 
for all $m\gg 0$, 
i.e., there exists $N\ge 1$ $($depending on $a$ and $b$$)$ 
such that $b a^m  \in M$ $($resp.,  $a^m b \in M$$)$ for all $m\ge N$.
\end{defi}   

A $R$-subspace $M$ of $\cA$ is said to be a {\it pre-two-sided Mathieu subspace} of $\cA$ if it is both left and right Mathieu subspace of $\cA$. Note that the {\it pre-two-sided} Mathieu subspaces were called  {\it two-sided} Mathieu subspace or {\it Mathieu subspaces} in \cite{GIC}. 
The change of the name here is due to the following family 
of two-sided Mathieu subspaces, which were not discussed 
in \cite{GIC} but are more entitled to be called 
(two-sided) Mathieu subspaces.  

\begin{defi}\label{Def-MS4}
A $R$-subspace $M$ of a  $R$-algebra $\cA$ is said to be 
a two-sided Mathieu subspace, or simply a Mathieu subspace,  
of $\cA$ if the following property holds: 
let $a\in \cA$ such that $a^m\in M$ for all $m\ge 1$. 
Then for any $b, c\in \cA$, we have 
$ b a^m c \in M$ for all $m\gg 0$, i.e., there exists $N\ge 1$ 
$($depending on $a$, $b$ and $c$$)$ such that 
$b a^m c \in M$ for all $m\ge N$.
\end{defi}   

Three remarks are as follows. First, all the algebras $\cA$  
involved in this paper are assumed to unital. 
For these algebras, it is easy to see that 
every (two-sided) Mathieu subspace is 
a pre-two-sided 
(and hence, also one-sided) Mathieu subspace.     

Second, from Definitions \ref{Def-MS} and \ref{Def-MS4} 
it is also easy to see that every left (resp., right) ideal 
of $\cA$ is a left (resp., right) Mathieu subspace of $\cA$, and 
every (two-sided) ideal of $\cA$ is a (two-sided) Mathieu subspace 
and hence, also a pre-two-sided Mathieu subspace of $\cA$. 
But the converse is not true 
(see \cite{DK}, \cite{GIC}, \cite{EWZ1}, \cite{FPYZ}, \cite{WZ} 
for some examples of Mathieu subspaces which 
are not ideals). Therefore, the notion of Mathieu subspaces 
can be viewed as a generalization of the notion of ideals.   

Third, just like the notion of ideals which has 
a generalization for modules of algebras,  
namely, the notion of submodules, 
the notion of Mathieu subspaces can also be generalized 
to modules of associative algebras. For more discussions 
in this direction, see \cite{GMS}. 

%Nevertheless, almost nothing is known in the literature 
%about this seemingly very mysterious concept.

The introduction of the notion of Mathieu subspaces in \cite{GIC} 
was mainly motivated by the studies of {\it the Mathieu conjecture}
 \cite{Ma}, {\it the vanishing conjecture} \cite{HNP}, \cite{GVC}, 
\cite{AGVC} and more recently, {\it the image conjecture} \cite{IC}, 
and also the well-known {\it Jacobian conjecture} \cite{K} 
(see also \cite{BCW}, \cite{E}). Actually, both 
{\it the Mathieu conjecture} and 
{\it the image conjecture} imply 
{\it the Jacobian conjecture}, and 
both are (open) problems on whether or not certain subspaces of some algebras are Mathieu subspaces 
(see \cite{Ma}, \cite{IC} and \cite{GIC} 
for more detailed discussions). The notion was named  after Olivier Mathieu due to his conjecture 
mentioned above.

There are also several other open problems and conjectures 
that are directly or indirectly related with Mathieu subspaces. 
For example, {\it the Dixmier conjecture} \cite{D} as shown first by  Y. Tsuchimoto \cite{Ts} in $2005$, and later by A. Belov and M. Kontsevich
\cite{BK} and P. K. Adjamagbo and A. van den Essen \cite{AE} in $2007$ 
is actually equivalent to {\it the Jacobian conjecture}; 
and {\it the vanishing conjecture} \cite{HNP}, \cite{GVC} on differential operators with constant coefficients, which now becomes a special case 
of {\it the image conjecture}, also implies {\it the Jacobian conjecture}.    
 
Furthermore, it has also been proposed in 
Conjecture $3.2$ in \cite{GIC}   
that the subspace of polynomials in $n\ge 1$ variables 
with complex coefficients whose integrals over a fixed open subset of $\bR^n$ with a positive measure are equal to zero  
should be a Mathieu subspace of the polynomial 
algebra in $n$ variables over $\bC$. 
In particular, by choosing some open subsets 
of $\bR^n$ and positive measures properly,  
this conjecture is equivalent to saying that 
every family of classical orthogonal polynomials 
(see \cite{Sz},  \cite{C}, \cite{DX}) in one or more variables 
with positive degrees should also span a co-dimension one 
Mathieu subspaces of the polynomial algebra 
(see Conjecture $3.5$ and the related discussions in \cite{GIC}). 
For some recent developments on the latter conjecture, 
see \cite{EWZ2}, \cite{FPYZ} and \cite{EZ}. For a recent 
survey on {\it the image conjecture} and its relations 
with {\it the vanishing conjecture}, {\it the Jacobian conjecture} and also the conjectures mentioned above, see \cite{E2}. 

Surprisingly, the conjecture on integrals of polynomials mentioned above is also related with the so-called {\it polynomial moment problem} proposed by M. Briskin, J.-P. Francoise and Y. Yomdin in the series of papers \cite{BFY1}-\cite{BFY5}, which was mainly motivated by the center problem for the complex Abel equation. For some recent studies on the {\it polynomial moment problem} in one or more variables, see \cite{PM}, \cite{Pa}, \cite{GIC} and \cite{FPYZ}.

Currently, it is also under investigations by the author  
and some of his colleagues whether or not images of all locally nilpotent derivations, locally finite derivations and divergence-zero derivations 
of polynomial algebras over fields of characteristic zero 
are Mathieu subspaces of the polynomial algebras. 
For example, it has been shown recently in \cite{EWZ1} that 
this is indeed the case for all locally finite derivations  
of polynomial algebras in two variables. It has also been 
shown in \cite{EWZ1} that for the two-variable case  
the same problem for the divergence-zero derivations 
having $1$ in the image is actually equivalent to 
the two-dimensional {\it Jacobian problem}. 
Furthermore, some Mathieu subspaces of 
the group algebras of finite groups have 
also been studied recently in \cite{WZ}. 

Due to their connections with the various open problems or conjectures mentioned above, especially their connections with {\it the Jacobian conjecture} and {\it the Dixmier conjecture}, the seemingly familiar but still very mysterious Mathieu subspaces deserve much more attentions from mathematicians. It is important and also necessary to study Mathieu subspaces in a separate and abstract setting.

\subsection{Contents and Arrangements} 
Before we proceed, one remark is in order. Even though 
most of the results on Mathieu spaces in this paper 
are stated and proved for all the four types 
({\it left}, {\it right}, {\it pre-two-sided} and {\it two-sided}) 
of Mathieu subspaces, for simplicity, 
in this subsection we only discuss the results    
for the {\it two-sided} case, i.e., 
only for Mathieu subspaces.  

In this paper, we first study some properties of  
the {\it radicals} of arbitrary subspaces and 
Mathieu subspaces of (associative) algebras, 
where for any $R$-subspace $V$ of a  $R$-algebras 
$\cA$, the {\it radical} of $V$, 
denoted by $\sqrt V$ or $\rad(V)$, 
is defined to be the set of the elements $a\in \cA$ such that 
$a^m\in V$ when $m\gg 0$. We then prove some properties and 
characterizations for the Mathieu subspaces with algebraic 
radicals for algebras over fields. 

One crucial result derived in this paper 
(see Theorem \ref{Char4AlgElt}) is that 
when the base ring $R$ is a field $K$, 
for algebraic elements $a\in \cA$, the positive integers   
$N$ in Definitions \ref{Def-MS} and \ref{Def-MS4} 
actually can be chosen in a way that does not depend 
on the element $b, c\in \cA$. Another crucial result 
for $K$-algebras $\cA$ is Theorem \ref{CharByIdem} 
which gives a characterization for Mathieu subspaces $V$ 
with algebraic radicals in terms of the 
idempotents contained in $V$.   
Consequently, for algebraic $K$-algebras, 
the Mathieu subspaces have an equivalent 
formulation that is much more similar to 
the definition of ideals 
(see Remark \ref{Def-MS3}).

By using some results derived in this paper, 
we also give characterizations or classifications 
for {\it strongly simple algebras} 
(see Definition \ref{StrSimAlg}) 
over arbitrary commutative rings, 
and for {\it quasi-stable algebras} 
(see Definition \ref{q-StaAlg}) 
over arbitrary fields 
(see Theorems \ref{No-MS-Thm1}, 
Proposition \ref{SpecialDomains} and 
Theorem \ref{Class-Qstable}). 
Furthermore, the co-dimension 
one Mathieu subspaces and the minimal non-trivial  Mathieu subspaces 
of all types are also completely classified for 
(finite dimensional) matrix algebras over fields 
(see Theorem \ref{Co-D1} and Proposition \ref{ClsMin}).  

Considering the length of this paper, below 
we give a more detailed description 
for the arrangements of the paper. 

In Section \ref{S2}, we first fix some notations and conventions 
that will be used throughout this paper. We then study certain properties of the radicals of Mathieu subspaces or arbitrary $R$-subspaces of $\cA$. A formally stronger but equivalent definition of Mathieu subspaces is also given in Proposition \ref{MS-Def2}.

In Section \ref{S3}, we study the algebraic elements of the radicals of  arbitrary subspaces $V$ and Mathieu subspaces $M$ of 
$K$-algebras $\cA$. The main results of this section are 
Theorems \ref{Idem-Thm1},  \ref{Alg-Ideal-Thm} and \ref{Char4AlgElt}. 
Theorem \ref{Idem-Thm1} says that $\sqrt V$ has no non-trivial idempotents of $\cA$ iff all algebraic elements of $\sqrt V$ are either nilpotent or invertible. 
Theorem \ref{Char4AlgElt} gives a characterization for  algebraic elements in the radicals of Mathieu subspaces 
$M$ of $\cA$, namely, for each algebraic $a\in \cA$, 
$a\in\sqrt M$ iff the principal ideal $(a^N) \subseteq M$ 
for some $N\ge 1$. Under the condition that $a^m\in M$ for 
all $m\ge 1$, Theorem \ref{Alg-Ideal-Thm} says that one can 
actually choose the integer $N$ above to be the multiplicity of $0\in K$ 
as a root of the minimal polynomial of the algebraic 
element $a\in \cA$.

In Section \ref{S4}, we use the results derived in 
Sections \ref{S2} and \ref{S3} to study various properties 
of Mathieu subspaces $M$ with algebraic radicals. For convenience, 
for any $K$-algebra $\cA$, we denote by $\gca$  
(resp.,  $\cE(\cA)$) the set of $K$-subspaces 
(resp.,  Mathieu subspaces) $V$ of $\cA$ such that 
$\sqrt V$ is algebraic over $K$, 
i.e., every element of $\sqrt V$ 
is algebraic over $K$.

In Subsection \ref{S4.1}, we give a characterization for 
Mathieu subspaces $V\in \cE(\cA)$ in terms of idempotents of $\cA$ 
(see Theorem \ref{CharByIdem}). Namely, 
a  $K$-subspace $V\in \gca$ is a Mathieu subspace of $\cA$ 
iff it contains the ideals of $\cA$ generated by 
the idempotents contained in $V$. In particular, 
the Mathieu subspaces of simple algebraic $K$-algebras 
can be characterized as $K$-subspaces of $\cA$ which do not 
contain any nonzero idempotents 
(see Proposition \ref{CharByIdem-cor2}).
Furthermore, the one-dimensional Mathieu subspaces of all 
$K$-algebras have been characterized in Proposition \ref{1D-MS}. 
This proposition will play some important roles in the later 
Sections \ref{S5}-\ref{S7}. 

In Subsection \ref{S4.2}, we study the relations between the 
radical of $M\in \cE(\cA)$ and the radical of the  
maximum ideal $I_M$ contained 
in $M$. In Lemma \ref{sqM=sqI}, and 
more generally in Theorem \ref{sqM=sqI-Thm}, 
we show that these two radicals actually coincide 
with each other.  In Theorem \ref{Comm-MSs}, we show that 
when $\cA$ is commutative, a  $K$-subspace $V\in \gca$ is 
a Mathieu subspace of $\cA$ iff its radical $\sqrt V$ 
is an ideal of $\cA$.

In Subsection \ref{S4.3}, we first show in Proposition 
\ref{Interscn} that the intersection of any family of Mathieu 
subspaces in $\cE(\cA)$ is still a Mathieu subspace of $\cA$. 
We then show in Proposition \ref{Union} that 
the union of any ascending sequence 
of Mathieu subspaces is also a Mathieu subspace of $\cA$ 
provided that the radical of the union is algebraic over $K$.   
Combining Propositions  \ref{Interscn} and  \ref{Union} 
with Zorn's lemma, we get existences of maximal or minimal  
elements in certain collections of Mathieu subspaces 
of algebraic $K$-algebras 
(see Proposition \ref{MaxMinMSs-Propo}, 
Theorem \ref{MaxMSs-Thm} and   
Corollary \ref{MaxMSs-Corol}).

In Section \ref{S5}, we show in Theorem \ref{Co-D1} that the only possible co-dimension one Mathieu subspace in the matrix algebra $M_n(K)$ $(n\ge 1)$ over a field $K$ is the subspace $H$ of the trace-zero matrices. More  precisely, if $char.\,K=p\le n$, $M_n(K)$ has no co-dimension 
one Mathieu subspace; and if $char.\,K=0$ 
or $char.\,K=p>n$, $H$ is the only co-dimension one 
Mathieu subspace (of any type) of $M_n(K)$. 
In Proposition \ref{ClsMin}) we show that the set of the nonzero 
minimal Mathieu subspaces is the same as the set of all dimension one $K$-subspaces of $M_n(K)$, 
which are not spanned by idempotent matrices.  

In Section \ref{S6}, we study the so-called {\it strongly simple algebras} 
$\cA$ over arbitrary commutative rings $R$, i.e., the $R$-algebras $\cA$  
whose only Mathieu subspaces are $0$ and $\cA$ itself.  
Note that every strongly simple algebra is a simple algebra 
since any ideal of $\cA$ is a Mathieu subspace 
of $\cA$. Under the convenient assumption $R\subseteq \cA$, we first show in Theorem \ref{No-MS-Thm1} that if a  $R$-algebra $\cA$ is strongly simple, 
then the base ring $R$ must be an integral domain and $\cA \simeq K_R$ 
as $R$-algebras, where $K_R$ denotes the field 
of fractions of $R$. In particular, for any field $K$, 
there are no strongly simple $K$-algebras except $K$ itself. 

We then show in Lemma \ref{ValuationLemma} that for every integral domain 
$R$ such that $R\ne K_R$ and $K_R$ has a real-valued additive valuation 
$\nu: K_R\to \bR$ satisfying $\nu(r)\ge 0$ for all $r\in R$, there is no strongly simple $R$-algebras. Note that this is the case for all Krull domains and Noetherian domains which are not fields   
(see Proposition \ref{SpecialDomains}). 
Consequently, all (commutative or noncommutative)  
rings except the finite fields $\bZ_p$ (for all primes $p$)  
are strongly simple $\bZ$-algebras   
(see Proposition \ref{RingCase}).    

In Section \ref{S7}, we first introduce the notions of 
({\it quasi}-){\it stable algebras} in Definition \ref{q-StaAlg}. 
We show in Proposition \ref{q-Quasi/R} that every integral $R$-algebra $\cA$, 
all of whose elements are either invertible or nilpotent, is quasi-stable. Consequently, every left or right integral Artinian local $R$-algebra 
is quasi-stable (see Corollary \ref{ArtinLocal}). 

We then give a classification in Theorem \ref{Class-Qstable} 
for the {\it quasi-stable algebras} over fields $K$. 
More precisely, we show that a  $K$-algebra $\cA$ 
is {\it quasi-stable} iff either $\cA\simeq K\dot{+}K$ 
or $\cA$ is an algebraic local $K$-algebra.
Note that by Corollary \ref{Alg-3Equiv}, 
the latter holds iff $\cA$ is algebraic and every element 
of $\cA$ is either nilpotent or invertible iff
$\cA$ is algebraic and has no non-trivial idempotents.

The motivation of the study of quasi-stable algebras 
is given in Proposition \ref{MotivStable} and 
Corollary \ref{MotivStable-Corol}. 
An application of Theorem \ref{Class-Qstable} 
via Corollary \ref{MotivStable-Corol} to commutative 
$K$-algebras is given in Corollary \ref{App2Noe}. Finally, 
for the completeness and also for the purpose of comparison, 
we also classify in Proposition \ref{ClassiStable} the 
{\it stable} $K$-algebras, i.e., the $K$-algebras $\cA$ 
such that every $K$-subspace $V\subset \cA$  
with $1\not \in V$ is an ideal of $\cA$.

\renewcommand{\theequation}{\thesection.\arabic{equation}}
\renewcommand{\therema}{\thesection.\arabic{rema}}
\setcounter{equation}{0}
\setcounter{rema}{0}

\section{\bf Mathieu Subspaces and Their Radicals}\label{S2}

In this section, we study some general properties 
of Mathieu subspaces and the radicals of subspaces of 
associative algebras. Most of the results derived 
in this section will be needed in the later sections. 

First, let's fix the following conventions 
and notations that will be used throughout 
this paper. 

Unless stated otherwise, $R$ and $K$ 
always stand for an arbitrary commutative 
ring and an arbitrary field, respectively. 
$\cA$ stands for an arbitrary    
associative (but not necessarily commutative) 
algebra over $R$ or $K$. Although most of the results 
in this paper also hold for non-unital algebras $\cA$, 
for convenience we assume that all rings and 
algebras in this paper have the identity elements 
which will be uniformly denoted by $1$, when no confusions 
occur. All algebra homomorphisms are assumed to 
preserve the identity elements. 
The ring or algebra with a single element $0$   
will be excluded in this paper.

Moreover, the following terminologies and notations 
for $R$-algebras $\cA$ will also be in force 
throughout this paper.

\begin{enumerate}
\item[1)] The sets of units or invertible elements 
of $R$ and $\cA$ will be denoted by $R^\times$ 
and $\cA^\times$, respectively.  

\item[2)] A $R$-subspace $V$ of $\cA$ is said to be {\it proper} 
if $V\ne \cA$, and {\it non-trivial} if $V\ne 0$ or $\cA$.

\item[3)] An element $a\in \cA$ is said to be an {\it idempotent} if $a^2=a$, 
and a {\it quasi-idempotent} if $a^2=r a$ 
for some $r\in R^\times$. An idempotent  
$a\in \cA$ is said to be {\it non-trivial} 
if $a\ne 0$ or $1\in \cA$.  

\item[4)] For any subset $S$ of a $R$-algebra $\cA$, 
we say $S$ is {\it integral} or {\it algebraic} 
(when $R$ is a field) over $R$ if every element 
$a\in S$ is integral over $R$ (i.e., $a$ 
is a root of a monic polynomial with 
coefficients in $R$). 

%\item We say that a  $R$-submodule or $R$-subspace 
%$V$ of $\cA$ is {\it non-trivial } if $V\ne \{0\}$ or $\cA$.

\item[5)] For any subset $S\subseteq \cA$, we define the 
{\it radical} of $S$, denoted by $\sqrt{S}$ or $\rad(S)$, 
to be the set of all the elements $a\in \cA$ such that
$a^m\in S$ when $m\gg 0$. The subset of the elements 
in the radical $\sqrt S$ which are integral over $R$ 
will be denoted by $\arad(S)$. 

\item[6)] The radical $\sqrt 0$ of the zero ideal will also be denoted by $\nil(\cA)$. Note that when $\cA$ is commutative, $\nil(\cA)$ is the {\it nilradical} of $\cA$. 

\item[7)] Let $\cA$ and $\cB$ be $R$-algebras.
We denote by $\cA\dot{+}\cB$ the $R$-algebra 
with the base $R$-space  $\cA \times \cB$ and the algebra product defined componentwise. 
\end{enumerate}

Note that for both Mathieu subspaces and ideals, 
we have several different cases:  
{\it left, right} and {\it $($pre-$)$two-sided}. 
Very often, it is necessary and important to treat all 
these cases. For simplicity, we introduce the short  
terminology {\it $\vartheta$-Mathieu subspaces} for 
Mathieu subspaces, where $\vartheta$ stands for 
{\it left, right}, {\it pre-two-sided}, or 
{\it two-sided}. Similarly, we introduce the terminology 
{\it $\vartheta$-ideals} for ideals, except for the 
specification  $\vartheta=${\it ``pre-two-sided\hspace{.1mm}"}, 
we also set {\it $\vartheta$-ideals} to mean 
{\it two-sided} ideals.    

In other words, the reader should read the letter $\vartheta$ as 
an index or a variable with four possible choices or 
``values".  However, to avoid repeating the phrase 
``{\it for every specification of $\vartheta$}" or {\it ``for every $\vt$"} 
infinitely many times, we will simply leave $\vartheta$ unspecified 
for the statements or propositions which hold 
for all the four specifications of 
$\vartheta$.  

Note that with the short terminologies 
fixed above, we immediately have the 
implication: {\it any $\vt$-ideal of $\cA$ is a $\vt$-Mathieu subspace of $\cA$}, which by the 
convention fixed above actually means 
four implications (corresponding to the four  
specifications of $\vt$).     

%For example, 
%\begin{enumerate}
%  \item[$i)$] 
%by saying ``{\it every $\vartheta$-ideal of $\cA$ is a 
%\vms of $\cA$}", we mean that all the four statements 
%obtained by replacing $\vartheta$ by  
%{\it left, right}, {\it pre-two-sided} and 
%{\it two-sided}, respectively, hold; 
% 
%\item[$ii)$] by saying  
%``{\it $\cA$ has at least one maximal non-trivial   
%$\vartheta$-Mathieu subspace}", we mean that 
%``{\it $\cA$ has at least one proper maximal 
%left Mathieu subspace, one proper maximal right Mathieu subspace, 
%one proper maximal pre-two-sided Mathieu subspace, and also one 
%proper maximal $($two-sided$)$ Mathieu subspace, 
%which are not necessarily the same as one another.}"
%\end{enumerate}
%
%Moreover, for each each element $S$ of $\cA$, 
%we denote by $(S)_\vartheta$ the $\vartheta$-ideal 
%of $\cA$ generated by elements of $S$.
%When $S$ consists of a single element 
%$a\in \cA$, we simply write the 
%$\vartheta$-ideal $(S)_\vartheta$ as $(a)_\vartheta$.
%However, for (two-sided) ideals, the commonly used notations 
%$(S)$ and $(a)$ will also be freely used. For instance, 
%$(a)=(a)_\vartheta$ with 
%$\vartheta=${\it two-sided}. \\

Finally, we fix the following notations. 

For any $a\in \cA$ and any $\vartheta\neq${\it ``pre-two-sided\hspace{.1mm}"}, 
we let $(a)_\vartheta$ denote the $\vartheta$-ideal 
of $\cA$ generated by $a$.  
For the case $\vartheta=${\it ``pre-two-sided"}, 
we set $(a)_\vartheta\!:=aA+Aa$, i.e., the sum of the left ideal 
and the right ideal generated by $a$. 
Moreover, for the two-sided case, the commonly used notation 
$(a)$ will also be freely used, i.e., $(a)=(a)_\vartheta$ 
with $\vartheta=$$\mbox{{\it ``two-sided\hspace{.1mm}"}}$. \\

Now let's start with the following formally stronger 
but equivalent definition of $\vartheta$-Mathieu subspaces, which says 
that the condition ``{\it $a^m\in M$ for all $m\ge 1$}" in 
Definitions \ref{Def-MS} and \ref{Def-MS4} may be replaced 
by the condition ``{\it $a\in \sqrt M\,$}". 

\begin{propo}\label{MS-Def2}
Let $\cA$ be a  $R$-algebra and $M$ a  $R$-subspace 
of $\cA$. Then $M$ is a $\vartheta$-Mathieu subspace 
of $\cA$ iff the following property holds: for any 
$a\in \sqrt M$ and $b, c\in \cA$, we have
\begin{enumerate}
  \item [$i)$]  $b a^m \in M$ when $m\gg 0$, if $\vt=${\it ``left"};

  \item [$ii)$]  $a^m c \in M$ when $m\gg 0$, if $\vt=${\it ``right"};
  
  \item [$iii)$] $b a^m, \, a^m c \in M$ when $m\gg 0$, 
             if $\vt=${\it ``pre-two-sided\hspace{.1mm}"};

  \item [$iv)$]  $b a^m c \in M$ when $m\gg 0$, 
  if $\vt=${\it ``two-sided\hspace{.1mm}"}. 
\end{enumerate} 
\end{propo}

\pf The $(\Leftarrow)$ part is trivial. To show 
the $(\Rightarrow)$ part, note first that since 
$a\in \sqrt M$, there exists $N \in \bN$ such that 
$a^m\in M$ for all $m\ge N$. 
%We need show that, for any $b\in \cA$, 
%$a^m b \in M$ when $m\gg 0$. 
Set $x\!:=a^N$. Then $x^m=a^{Nm}\in M$ for all $m\ge 1$. 

Assume that $M$ is a (two-sided) Mathieu subspace of $\cA$. 
Then for any $b, c\in \cA$, by Definition \ref{Def-MS4}  
it is easy to see that 
for the (finitely many) elements $b a^r  \in \cA$ 
$(0\le r\le N-1)$, there exists 
$N_1\in \bN$ such that 
\begin{align*} 
b a^{Nm+r} c =(b a^r)x^m c \in M 
\end{align*}
for all $0\le r\le N-1$ and $m\ge N_1$.

From the equation above, it is easy to see 
that for all $k\ge NN_1$, we have $b a^k c \in M$.
Therefore, the theorem holds for (two-sided)  
Mathieu subspaces. 

By letting $c=1$ (resp., $b=1$) in the arguments above, 
we see that the theorem also holds for left 
(resp., right) Mathieu subspaces, whence 
the pre-two-sided case also follows.   
\epfv

Next, we use a similar argument as in the proof above 
to show the following lemma on the radicals
of $\vartheta$-Mathieu subspaces. 

\begin{lemma}\label{General-RadLemma}
Let $\cA$ be a  $R$-algebra and $S$ 
a subset of $\cA$. Then the following 
statements hold.
\begin{enumerate}
  \item [$i)$] $\sqrt S \subseteq \rad(\sqrt S)$.

  \item [$ii)$] Assume further that $S$ is a \vms of $\cA$. 
Then $\sqrt S = \rad(\sqrt S)$.  
\end{enumerate} 
\end{lemma}  

\pf $i)$ Let $a\in \sqrt S$. Then we have $a^m \in S$ when 
$m\gg 0$. Hence, for any $k\ge 1$, we also have 
$(a^k)^m=a^{km}\in S$ when $m\gg 0$, whence 
$a^k\in \sqrt S$. Therefore, $a\in \rad(\sqrt S)$  
and hence, the statement follows.

$ii)$ Let $a\in \rad(\sqrt S)$. 
Then $a^m\in \sqrt S$ when $m\gg 0$, 
i.e., there exists $N\ge 1$ 
such that $a^N \in  \sqrt S$. 
Since $S$ is a \vms of $\cA$, 
by Proposition \ref{MS-Def2} 
there exists $N_1\ge 1$ such that 
$a^{Nm+r}=(a^N)^m a^r \in S$ 
for all $0\le r\le N-1$ 
when $m\ge N_1$. From this fact 
it is easy to see that for all 
$k\ge NN_1$, we have $a^k\in S$. 
Therefore, $a\in \sqrt S$  
and hence, $\rad{(\sqrt S)} \subseteq \sqrt S$.
Then by $i)$, the equality in $ii)$ 
follows.
\epfv

Note that the statement $ii)$ in Lemma \ref{General-RadLemma} 
is parallel to the fact in commutative algebra 
that the radicals of ideals are radical. 
Of course, in general the radicals of \vmss of $\cA$  
are not closed under the addition or the product of 
the algebra $\cA$. But, as we will see later in 
Theorem \ref{Comm-MSs} and Corollary \ref{Comm-MSs-Corol}, 
for commutative $K$-algebras $\cA$, 
a  $K$-subspace $V\subseteq \cA$ with  
$\sqrt V$ algebraic over $K$ is a Mathieu subspace 
of $\cA$ iff its radical $\sqrt V$ is a radical 
ideal of $\cA$. 

Next, we give the following characterizations for  
$\vartheta$-ideals and $\vartheta$-Mathieu 
subspaces. 
Since Eqs.\,(\ref{CharsInRads-e4}), (\ref{CharsInRads-e3}) and
(\ref{CharsInRads-e7}) (below) 
obviously imply 
Eqs.\,(\ref{CharsInRads-e6}), (\ref{CharsInRads-e5}) and 
(\ref{CharsInRads-e8}), respectively, 
the characterizations provide  
a different point of view to see that 
the notion of \vmss is indeed 
a natural generalization of 
the notion of $\vartheta$-ideals.

\begin{lemma}\label{CharsInRads}
Let $V$ be a  $R$-subspace of a  $R$-algebra 
$\cA$. For each $b\in \cA$, we set 
\allowdisplaybreaks{
\begin{align}
(V:b)\!:&=\{a\in \cA\,|\, ab \in V\}, \label{CharsInRads-e1} \\ 
b^{-1}V\!:&=\{a\in \cA\,|\, ba \in V\}, \label{CharsInRads-e2}
\end{align} }
where $b^{-1}V$ is an abusing notation 
since $b$ might not be invertible in $\cA$. 

Then the following statements hold.
\begin{enumerate} 
\item[$i)$] $V$ is a left ideal of $\cA$ iff 
 for any $b\in \cA$, we have 
  \begin{align}\label{CharsInRads-e4}
V \subseteq b^{-1}V.
\end{align}

\item[$ii)$] $V$ is a left Mathieu subspace of $\cA$ iff 
 for any $b\in \cA$, we have 
  \begin{align}\label{CharsInRads-e6}
\sqrt{V} \subseteq  \sqrt{b^{-1}V}.
\end{align}

 \item[$iii)$] $V$ is a right ideal of $\cA$ iff 
for any $b\in \cA$, we have
  \begin{align}\label{CharsInRads-e3}
V\subseteq (V: b).
\end{align}

\item[$iv)$] $V$ is a right Mathieu subspace of 
$\cA$ iff for any $b\in \cA$, we have 
  \begin{align}\label{CharsInRads-e5}
\sqrt{V} \subseteq \sqrt{(V: b)}.
\end{align}

\item[$v)$] $V$ is a $($two-sided$)$ ideal of $\cA$ iff 
for any $b, c \in \cA$, we have
  \begin{align}\label{CharsInRads-e7}
V\subseteq b^{-1}(V: c).
\end{align}

\item[$vi)$] $V$ is a $($two-sided$)$ Mathieu subspace of 
$\cA$ iff for any $b, c\in \cA$, we have 
  \begin{align}\label{CharsInRads-e8}
\sqrt{V} \subseteq \sqrt{b^{-1}(V: c)}.
\end{align}
\end{enumerate}
\end{lemma}
\pf The proof of the lemma is very straightforward. Here we just give a proof for the statement $vi)$. The other statements can be proved similarly. 

$(\Leftarrow)$ Let $a\in \sqrt V$. Then by 
Eq.\,(\ref{CharsInRads-e8}), $a\in \sqrt{b^{-1}(V:c)}$, i.e., 
$a^m\in b^{-1}(V:c)$ when $m\gg 0$. Hence, by 
Eqs.\,(\ref{CharsInRads-e1}) and (\ref{CharsInRads-e2}), 
we have $ba^mc\in V$ when $m\gg 0$. It then follows from 
Proposition \ref{MS-Def2} that $V$ is 
a $($two-sided$)$ Mathieu subspace.

The $(\Rightarrow)$ part follows simply by reversing 
the arguments above.
\epfv

Next, we prove another lemma on 
the radicals of $R$-subspaces of $\cA$, 
which will be needed later in 
Subsection \ref{S4.2}.

\begin{lemma}\label{CuteLemma}
Let $\cA$ be a  $R$-algebra $($not necessarily commutative$)$ 
and $V$ a  $R$-subspace of $\cA$ such 
that $\sqrt V=\cA$. 
Then $V=\cA$.
\end{lemma}

\pf Assume otherwise and let $a\in \cA\backslash V$. 
Since $a\in \cA=\sqrt V$, we have $a^m\in V$ when 
$m\gg 0$. Since $a\not \in V$,  
there exists $k\ge 1$ such that 
$a^k\not \in V$ but $a^m\in V$ 
for all $m\ge k+1$. 

Set $b\!:=1+a^k$. Since $b\in \sqrt V(=\cA)$, 
there exists $N\ge 1$ 
such that $b^m\in V$ for all $m\ge N$. 
Note that for each $m\ge N$, we also have 
\begin{align}
b^m=(1+a^k)^m \equiv 1+ma^k \mod V.
\end{align}
Therefore, $1+ma^k\in V$ for all $m\ge N$. 
Consequently, we have 
$a^k=(1+(N+1)a^k)-(1+Na^k)\in V$, 
which is a contradiction.  
\epfv

Now let's recall the following simple but 
very useful property of $\vartheta$-Mathieu 
subspaces, which can be easily checked  
(or see Proposition 4.9 in \cite{GIC}). 

%which was first proved in 
%Proposition 4.9 in \cite{GIC} for all 
%the cases except the two-sided ones. 
%But, it is easy to see that the proof 
%given in \cite{GIC} also works for 
%two-sided Mathieu subspaces.   

\begin{propo}\label{Pull-Back}
Let $\cA$ and $\cB$ be $R$-algebras  
and $\phi:\cA\to \cB$ a  $R$-algebra homomorphism. 
Then for every \vms $M$ of $\cB$, 
$\phi^{-1}(M)$ is a \vms of $\cA$.
\end{propo}

One immediate consequence of the proposition above 
is the following corollary.

\begin{corol}\label{Restriction}
Let $\cB$ be a  $R$-algebra and $\cA$ 
a  $R$-subalgebra of $\cB$. Then 
for every \vms $M$ of $\cB$, 
$M\cap \cA$ is a \vms of $\cA$.
\end{corol}
\pf Apply Proposition \ref{Pull-Back} to 
the embedding $\iota :\cA \to \cB$ and 
note that $\iota^{-1}(M)=M\cap \cA$.
\epfv 

\begin{propo}\label{quotient-I}
Let $I$ be an ideal of $\cA$ and 
$M$ a  $R$-subspace of $\cA$. 
Assume that $I\subseteq M$. 
Then $M$ is a \vms of 
$\cA$ iff $M/I$ is a \vms  
of $\cA/I$.
\end{propo}

\pf $(\Leftarrow)$ Let $\pi: \cA\to \cA/I$ be 
the quotient map. Since $I\subseteq  M$, 
we have $\pi^{-1}(M/I)=M+I=M$. Applying 
Proposition \ref{Pull-Back} to the 
$R$-algebra homomorphism $\pi$, we see 
that $M$ is a \vms of $\cA$.

$(\Rightarrow)$ Let 
$\bar a, \bar b \in \cA/I$ 
such that $\bar a^m\in M/I$ 
for all $m\ge 1$. Let $a, b\in \cA$ such that 
$\pi(a)=\bar a$ and $\pi(b)=\bar b$. 
Then for all $m\ge 1$,  we have $a^m\in \pi^{-1}(M/I)=M$  
since $\pi(a^m)=\bar a^m\in M/I$. 

Now assume $\vartheta=${\it``left"},  i.e., $M$ is a 
left Mathieu subspace of $\cA$. Then we have $b a^m \in M$ when $m\gg 0$, 
whence $\bar b \bar a^m =\pi(b a^m)\in M/I$ 
when $m\gg 0$. Therefore, $M/I$ is a left  
Mathieu subspace of $\cA/I$. For the other 
specifications of $\vartheta$, the proofs 
are similar.  
\epfv

The following lemma is obvious but does provide  
a family of $\vt$-Mathieu subspaces.

\begin{lemma}\label{radical-Lemma}
Let $M$ be a  $R$-subspace of $\cA$ such that  
$\sqrt{M}\subseteq  \nil(\cA)$. Then every $R$-subspace 
$V\subseteq M$ is a $\vt$-Mathieu subspace of $\cA$.
\end{lemma}
%\pf Since $V\subseteq M$, we have  
%$\sqrt{V}\subseteq \sqrt{M}\subseteq \nil (\cA)$, 
%from which and Definition \ref{MS-Def} 
%the lemma follows. 
%\epfv

The following lemma will be crucial for our later arguments.

\begin{lemma}\label{CyclicLemma}
Let $a$ be a nonzero quasi-idempotent of $\cA$ 
and $V$ a  $R$-subspace of $\cA$. 
Then the following statements hold. 
\begin{enumerate}
\item[$i)$] $a$ is integral over $R$ but cannot be nilpotent. 
Moreover, $a$ is invertible iff $a$ is an invertible scalar of $\cA$, 
i.e., $a \in R^\times \cdot 1_\cA\subset \cA$.

  \item[$ii)$]  $a\in \sqrt V$ iff $a\in V$.

  \item[$iii)$]  Assume further that $a\in V$ and  
$V$ is a \vms of $\cA$. 
Then $(a)_\vartheta \subseteq V$.
\end{enumerate}
\end{lemma}

\pf Assume $a^2=r a$ for some $r\in R^\times$. 
Then it follows inductively that $a^m=r^{m-1} a$ 
for all $m\ge 1$, from which it is easy see 
that $ii)$ does hold.

To show $i)$, note first that $a$ is 
integral over $R$ since $a$ is a root of the monic polynomial 
$t^2-rt=0$, and $a$ cannot be nilpotent, for if $a^m=0$ 
for some $m\ge 2$, then $a=r^{1-m}a^m=0$, 
which is a contradiction. Furthermore, if 
$a \in \cA^\times$, then 
from the equation $a(a-r)=0$, 
we have $a=r\in R^\times$. Since every invertible scalar of $\cA$  
is a quasi-idempotent, we see that $i)$ follows.

To show $iii)$, note first that by $ii)$ $a\in \sqrt V$. 
If $V$ is a left Mathieu subspace of $\cA$,  
then for each $b\in \cA$, we have $r^{m-1}ba=ba^m\in V$ 
when $m\gg 0$. Since $r^{m-1}\in R^\times$ for 
all $m\ge 1$, we have $ba\in V$, whence 
$\cA a\subseteq V$. 

The {\it right} and {\it two-sided} cases can be proved similarly. 
The {\it pre-two-sided} case follows directly from the  
{\it left} and {\it right} cases. 
\epfv

Applying Lemma \ref{CyclicLemma}, $iii)$ to the identity element 
$1\in \cA$, we immediately get the following corollary,  
which was first noticed in \cite{GIC}.

\begin{corol}\label{OneLemma}
For any \vms $M$ of $\cA$ with  
$1\in M$, we have $M=\cA$. 

Equivalently, any proper $R$-subspace 
$V\subset \cA$ with $1 \in V$ cannot be a 
$\vartheta$-Mathieu subspace of $\cA$.
\end{corol}

\renewcommand{\theequation}{\thesection.\arabic{equation}}
\renewcommand{\therema}{\thesection.\arabic{rema}}
\setcounter{equation}{0}
\setcounter{rema}{0}

\section{\bf Algebraic Elements in the Radicals  
of Arbitrary Subspaces}\label{S3}

In this section, we study some properties 
of integral or algebraic elements in the radicals 
of arbitrary subspaces 
or $\vartheta$-Mathieu subspaces of associative 
algebras $\cA$ over a commutative ring 
$R$ or a field $K$. 

Recall that for any subset $S\subseteq \cA$, we have 
let $\arad(S)$ to denote the subset of integral 
or algebraic (if the base ring $R$ is a field) elements 
in the radical $\sqrt S$ (or $\rad(S)$) 
of $S$. 

\begin{lemma}\label{Inv-Alg-L1-New}
Let $\cA$ be a  $R$-algebra and $V$ an arbitrary 
$R$-subspace of $\cA$. Assume that there exists $a\in \sqrt V$ 
such that $a$ is invertible and $a^{-1}$ is integral over $R$. 
Then $1\in V$.  
\end{lemma}

\pf Note first that by replacing $a$ by a positive power of $a$ 
if necessary, we may assume that $a^m\in V$ 
for all $m\ge 1$. 

Let $f(t)$ be a monic polynomial with coefficients 
in $R$ such that $f(a^{-1})=0$. 
Write $f(t)=t^d-\sum_{k=0}^{d-1} r_k t^k$ for 
some $d\ge 1$ and $r_k\in R$ $(0\le k\le d-1)$. 
Then we have 
\begin{align*}
a^{-d}-r_{d-1}a^{1-d}-r_{d-2}a^{2-d}-\cdots - r_1 a^{-1}-r_0=0.
\end{align*}

Multiplying $a^d$ to the equation above,  we get  
\begin{align*}
1=r_{d-1}a+ r_{d-2}a^2 \cdots + r_1 a^{d-1}+r_0 a^d.
\end{align*}

Since $a^m\in V$ for all $m\ge 1$, it follows from the equation 
above that $1\in V$. 
\epfv   

\begin{propo}\label{Inv-Alg-L1}
Let $K$ be a field, $\cA$ a  $K$-algebra and $V$ a $K$-subspace of $\cA$. 
Then $\arad(V)\cap\cA^\times \ne \emptyset$ 
iff $1\in V$. 
\end{propo}

\pf $(\Leftarrow)$ Since $1\in V$, then $1\in \sqrt V$. 
Since $1\in \cA$ is invertible and algebraic 
over $K$, we have $1\in \arad(V)\cap\cA^\times$, whence 
$\arad(V)\cap\cA^\times \ne \emptyset$. 

$(\Rightarrow)$ Let $a\in \arad(V)\cap\cA^\times$. Then  
$a\in \sqrt V$ and is invertible and algebraic over $K$. 
Since the base ring is a field $K$, it is easy to see that 
$a^{-1}$ is also algebraic over $K$. Then by 
Lemma \ref{Inv-Alg-L1-New}, we have $1\in V$.
\epfv

In order to get more results on 
algebraic elements in the radicals of 
$K$-subspaces of $K$-algebras, we need the following lemma 
on polynomials $f(t)$ in one variable $t$ over a field $K$.

\begin{lemma}\label{CrucialLemma}
Let $f(t)=t^k h(t)$ for some $k\ge 0$ and 
$h(t)\in K[t]$ such that $h(0)\ne 0$.
Then there exists a polynomial $p(t)\in K[t]$ 
such that the following equations hold: 
\allowdisplaybreaks{
\begin{align}
p(t)& \equiv 0 \qquad \mod (t^k), \label{CrucialLemma-e1} \\
p^2(t)& \equiv p(t) \quad \mod (f(t)), \label{CrucialLemma-e2}\\
t^k & \equiv t^k p(t)\, \mod (f(t)). \label{CrucialLemma-e3}
\end{align}
Furthermore, if $k\ge 1$ and $\deg h\ge 1$, we have
\begin{align}
p(t) \not \equiv 0, 1 \mod (f(t)). \label{CrucialLemma-e4} 
\end{align} }
\end{lemma}

\pf First, if $k=0$, we choose $p(t)=1$. Then 
it is easy to see that 
Eqs.\,(\ref{CrucialLemma-e1})--(\ref{CrucialLemma-e3}) in the 
lemma hold in this case.

Assume $k\ge 1$. Since $h(0)\ne 0$, the polynomials $t^k$ and 
$h(t)$ are co-prime. Therefore, there exist $u(t), v(t)\in K[t]$ such that 
\begin{align}
1= t^k u(t)+h(t)v(t). \label{CrucialLemma-pe1} 
\end{align}

Let $p(t)\!:=t^k u(t)$. Then Eq.(\ref{CrucialLemma-e1}) 
follows immediately. Furthermore, from Eq.(\ref{CrucialLemma-pe1}) 
we have 
\begin{align}
p(t)=1-h(t)v(t). \label{CrucialLemma-pe2} 
\end{align}

Multiplying $p(t)$ and $t^k$ to the both sides 
of the equation above, respectively, we get 
\begin{align}
p^2(t)&=p(t)-p(t)h(t)v(t)=p(t)- t^k u(t) h(t) v(t) \label{CrucialLemma-pe3} \\
&=p(t)-u(t)v(t)f(t), \nno \\
t^k p(t)&=t^k-t^kh(t)v(t)=t^k-f(t)v(t). \label{CrucialLemma-pe4} 
\end{align}

Then Eqs.\,(\ref{CrucialLemma-e2}) and (\ref{CrucialLemma-e3}) 
follow immediately from 
Eqs\,.(\ref{CrucialLemma-pe3}) and (\ref{CrucialLemma-pe4}), 
respectively.

Finally, we prove Eq.\,(\ref{CrucialLemma-e4}) as follows. 

Assume $p(t)\equiv 0 \,\,{\rm mod}\, (f(t))$. Then $f(t) \,|\, p(t)$, whence  
$h(t)\, | \, p(t)$. However, by Eq.\,(\ref{CrucialLemma-pe2}), 
we have $h(t)\, |\, 1$, which contradicts the condition 
$\deg h\ge 1$.

Assume $p(t)\equiv 1 \,\,{\rm mod}\, (f(t))$. Then we have $f(t)\,|\,(p(t)-1)$, 
whence  $t^k\, | \, (p(t)-1)$. By Eq.\,(\ref{CrucialLemma-pe2}), 
we also have $t^k\, |\, h(t)v(t)$. Hence $t^k\,|\, v(t)$ since 
$h(0)\ne 0$. Then by Eq.\,(\ref{CrucialLemma-pe1}), 
we get $t^k\,|\, 1$, which contradicts 
the condition $k\ge 1$.
\epfv

%\begin{exam}
%Let $f(t)=t(t^2-1)$. ....
%\end{exam}

\begin{propo}\label{Propo4p(a)}
Let $V$ be a  $K$-subspace of $\cA$ and 
$a$ an algebraic element of $\cA$, which is not nilpotent 
nor invertible. Denote by $k(\ge 1)$ the multiplicity 
of $0\in K$ as a root of the minimal polynomial  
of $a$ over $K$. Assume further that 
$a^m\in V$ for all $m\ge 1$. 
Then there exists $p(t)\in t^k K[t]$ such that 
the following three statements hold:
\begin{enumerate}
\item[$i)$] $p(a)\in V$;
  
\item[$ii)$] $p(a)$ is a non-trivial idempotent of $\cA$;

\item[$iii)$] $a^k=a^kp(a)$.  
\end{enumerate}  
\end{propo}

\pf Let $f(t)$ be the minimal polynomial  
of $a$ over $K$ and write it as 
$f(t)=t^k h(t)$ for some $h(t)\in K[t]$ 
such that $h(0)\ne 0$. 
Since $a$ is not nilpotent, we have $\deg h\ge 1$. 
Since $a$ is not invertible, we have $f(0)=0$, 
which means $k\ge 1$ as already indicated 
in the theorem.

Now apply Lemma \ref{CrucialLemma} to 
the polynomial $f(t)$ and let $p(t)$ be 
as in the same lemma. Then by Eq.\,(\ref{CrucialLemma-e1}), 
$p(t)\in t^kK[t]$, and by 
Eqs.\,(\ref{CrucialLemma-e2})-(\ref{CrucialLemma-e4}), 
$p(a)$ satisfies $ii)$ and $iii)$. 
To show $i)$, note that 
$p(t)\in t^k K[t]$ with $k\ge 1$. 
So $p(a)$ is a linear combination of some powers 
$a^m$'s over $K$ with $m\ge 1$. 
Since by our assumption  
$a^m\in V$ for all $m\ge 1$, 
we have $p(a)\in V$. 
\epfv

\begin{theo}\label{Idem-Thm1}
Let $\cA$ be a  $K$-algebra and $V$ a $K$-subspace of $\cA$. 
Then the following two statements 
are equivalent.

$1)$ every element of $\arad (V)$ is  
either nilpotent or invertible.  

$2)$ $V$ contains no non-trivial idempotents. 
\end{theo}

\pf $1)\Rightarrow 2)$: Assume that $V$ contains 
a non-trivial idempotent $e$. Then by Lemma \ref{CyclicLemma}, 
$i)$ and $ii)$,  
we know that $e\in \arad (V)$ and 
$e$ is not nilpotent nor invertible, 
which contradicts $1)$.  

$2)\Rightarrow 1)$: Assume that there exists 
$a\in \arad (V)$ which is not nilpotent 
nor invertible. Note that for each $m\ge 1$, 
$a^m$ is also algebraic over $K$ 
and is not nilpotent nor invertible. 
Since $a^m\in V$ when $m\gg 0$, replacing 
$a$ by a power of $a$ if necessary, 
we may further assume that $a^m\in V$ 
for all $m\ge 1$. Then 
by Proposition \ref{Propo4p(a)}, 
we get a non-trivial idempotent 
$p(a)\in V$, which is a contradiction.
\epfv

Applying the theorem above to $V=\cA$ and noting that 
$\sqrt{\cA}=\cA$, we immediately have the following 
corollary.

\begin{corol}\label{2Equiv}
For every $K$-algebra $\cA$, the following two statements are 
equivalent.

$1)$ every algebraic element of $\cA$ is either 
nilpotent or invertible. 

$2)$ $\cA$ has no non-trivial idempotents. 
\end{corol}

The following lemma and the corollary followed 
provide more understandings on the equivalent conditions 
in the corollary above. 

\begin{lemma}\label{3Equiv}
Let $\cA$ be a $R$-algebra. Then for the following three statements:  
\begin{enumerate}
  \item[$1)$] every element of $\cA$ is either nilpotent or invertible;
  \item[$2)$] $\cA$ is a local $R$-algebra;
  \item[$3)$] $\cA$ has no non-trivial idempotent,
  \end{enumerate}      
we have $1)\Rightarrow 2) \Rightarrow 3)$.
\end{lemma} 
\pf $1)\Rightarrow 2)$ is well-known, e.g., see Corollary $a$, 
p.\,74 in \cite{P}.
To show  $2)\Rightarrow 3)$, let ${\bf J}(\cA)$ be 
the Jacobson radical of $\cA$. Then ${\bf J}(\cA)$ 
is also the unique maximal left ideal of $\cA$. 
Assume that $\cA$ has a non-trivial idempotent $e$. 
Then it is easy to check that $1-e$ is 
also a non-trivial idempotent. 
Furthermore, by Lemma \ref{CyclicLemma}, $i)$, 
both $e$ and $1-e$ are not invertible, whence  
the left ideals $e\cA$ and $(1-e)\cA$ 
are proper and hence, both are contained 
in ${\bf J}(\cA)$. In particular, 
both $e$ and $1-e$ are in ${\bf J}(\cA)$. 
But this implies  $1=e+(1-e)\in {\bf J}(\cA)$, 
which is a contradiction. 
\epfv

\begin{corol}\label{Alg-3Equiv}
For every algebraic $K$-algebra $\cA$, the three statements  
in Lemma \ref{3Equiv} are equivalent to one another.  
\end{corol} 
\pf Since $\cA$ is algebraic, we see by Corollary \ref{2Equiv} that the statements $1)$ and $3)$ in Lemma \ref{3Equiv} 
are actually equivalent to each other. With this observation 
the corollary follows immediately from Lemma \ref{3Equiv}.
\epfv

%\begin{lemma}
%Every algebraic Artinian local $K$-algebra $\cA$ 
%(not necessarily commutative) satisfies 
%the equivalent statements Corollary \ref{2Equiv}. 
%\end{lemma}
%\pf Let $\frak m$ be the unique maximum ideal of 

%% Good Example!!

%However, not all algebraic $K$-algebras satisfying 
%the equivalent statements in Corollary \ref{2Equiv} 
%are Artinian. For example, let $\cA=K[x_i\,|\, i\ge 1]$ 
%be the polynomial algebra over $K$ in infinitely 
%many commutative variables $x_i$ $(i\ge 1)$, and $I$ 
%the ideal of $\cA$ generated by $x_i^i$ $(i\ge 1)$.
%Set $\cB\!:=\cA/I$. Then it is easy to see that $\cB$ 
%is algebraic over $B$, and satisfies statement $1)$ in 
%Corollary \ref{2Equiv}. But $\cB$ is not Noetherian, 
%and hence not Artinian, for the maximum ideal of $\cB$, 
%which is generated by the images of $x_i$ 
%$(i\ge 1)$ in $\cB$, is not finitely generated.

Next, we derive the following theorem on 
algebraic elements of the radicals of 
$\vartheta$-Mathieu subspaces.

\begin{theo}\label{Alg-Ideal-Thm}
Let $\cA$ be a  $K$-algebra and 
$M$ a \vms of $\cA$. 
Let $a\in \cA$ such that $a$ is algebraic over $K$ 
and $a^m\in M$ for all $m\ge 1$. 
Denote by $k\ge 0$ the multiplicity of $0\in K$ 
as a root of the minimal polynomial $f(t)$ of $a$. 
Then $(a^k)_\vartheta \subseteq M$. In particular, 
for any $\vt\neq${\it ``pre-two-sided\hspace{.1mm}"},   
the $\vartheta$-ideal of $\cA$ generated 
by $a^k$ is contained in $M$. 
\end{theo}

\pf Assume first that $k=0$, i.e., $0$ is not a root of 
$f(t)$. Then $a$ is invertible, and by 
Proposition \ref{Inv-Alg-L1}, $1\in M$. By Corollary \ref{OneLemma}, 
we have $M=\cA$. Hence the theorem holds in this case. 
 
Assume that $k\ge 1$. Then $a$ is not invertible. 
If $a$ is nilpotent, then $a^k=0$, whence  
the theorem holds trivially in this case. 
So assume that $a$ is not nilpotent nor invertible. Applying Proposition \ref{Propo4p(a)} to $a$ with $V=M$, and 
letting $p(a)$ be as in the same 
proposition, we see that $a^k=a^kp(a)$ and 
$p(a)$ is a non-trivial 
idempotent in $M$. 

Now, applying Lemma \ref{CyclicLemma}, $iii)$ to the idempotent 
$p(a)$ with $V=M$, we get $(p(a))_\vartheta \subseteq M$. 
Furthermore, since $a^k=a^kp(a)$, we also have
\begin{align}
(a^k)_\vartheta =(a^kp(a))_\vartheta 
\subseteq (p(a))_\vartheta \subseteq M.
\end{align}
Hence the theorem follows.
\epfv

One immediate consequence of Theorem \ref{Alg-Ideal-Thm}
is the following characterization of algebraic elements 
in the radicals of $\vartheta$-Mathieu subspaces.
 
\begin{theo}\label{Char4AlgElt}
Let $M$ be a \vms of a  $K$-algebra $\cA$ and 
$a$ an algebraic element of $\cA$. 
Then $a\in \sqrt{M}$ iff $(a^N)_\vartheta 
\subseteq M$ for some $N\ge 0$.
\end{theo}

\pf The $(\Leftarrow)$ part follows directly from the 
fact that for all $m\ge N$, 
$a^m \in (a^N)_\vartheta \subseteq M$.
The $(\Rightarrow)$ part can be proved 
as follows. 

Since $a\in \sqrt M$, we have that $a^m\in M$ when 
$m\gg 0$. In particular, there exists 
$n\ge 1$ such that $(a^n)^m=a^{nm}\in M$ for all $m\ge 1$.  
Applying Theorem \ref{Alg-Ideal-Thm} to the algebraic element 
$a^n\in \cA$, we have  
$(a^{nk})_\vartheta =\left((a^n)^k\right)_\vartheta \subseteq M$ 
for some $k\ge 0$, whence the theorem follows with $N=nk$. 
\epfv

\renewcommand{\theequation}{\thesection.\arabic{equation}}
\renewcommand{\therema}{\thesection.\arabic{rema}}
\setcounter{equation}{0}
\setcounter{rema}{0}

\section{\bf Mathieu Subspaces with Algebraic Radicals}
\label{S4}

Throughout this section, $K$ stands for an arbitrary  
field and $\cA$ an associative algebra 
over $K$. For convenience, we denote by $\gca$  
(resp., $\ovca$) the collection 
of all $K$-subspaces (resp., $\vartheta$-Mathieu subspaces) 
$V$ of $\cA$ such that $\sqrt V$ is algebraic over $K$. 
 
In this section we use the results derived in the previous sections to 
study some properties of \vmss  in $\ovca$. 
Note that all the results derived in this section apply 
under one of the conditions in the following 
easy-to-check lemma. 

\begin{lemma}\label{3Cases}
Let $V$ be a  $K$-subspace of $\cA$. Then 
$V\in \gca$  if one of the 
following four conditions holds:
\begin{enumerate}
  \item[a)]  $\cA$ is algebraic over $K$; 
  \item[b)]  $V$ is algebraic over $K$;
  \item[c)]  $\dim_K \cA < \infty$.
  \item[d)]  $\dim_K V < \infty$.
\end{enumerate}   
\end{lemma}

\subsection{Characterization of $M\in \ovca$ in Terms of Idempotents} \label{S4.1}

We start with the following characterization of \vmss 
in $\ovca$ in terms of idempotents of $\cA$.

\begin{theo}\label{CharByIdem}
Let $V\in \gca$. Then $V$ is a \vms of $\cA$ iff 
for any idempotent $e \in V$, we have $(e)_\vartheta \subseteq V$.
\end{theo}

\pf The $(\Rightarrow)$ part follows directly from 
Lemma \ref{CyclicLemma}, $iii)$. For the 
$(\Leftarrow)$ part, we here just give a proof for  
the two-sided case. The proofs for the other three 
cases are similar. 

Let $a, b, c\in \cA$ such that $a^m\in V$ for all $m\ge 1$. 
We need to show that $b a^m c\in V$ when $m\gg 0$.

Note first that since $a\in \sqrt V$ and $V\in \gca$, 
$a$ is algebraic over $K$. If $a$ is nilpotent, then 
$ba^mc=0\in V$ when $m\gg 0$. If $a$ is invertible, 
then $1\in V$ by Proposition \ref{Inv-Alg-L1}. 
Applying our assumption to the idempotent $1\in V$, 
we have $V=\cA$, whence $ba^mc\in V$ for all $m\ge 1$.

Finally, assume that $a$ is not nilpotent nor invertible. 
Apply Proposition \ref{Propo4p(a)} to $a$, and let $p(a)$ 
and $k\ge 1$ be as in the same proposition. 
Then $p(a)$ is an idempotent in $V$, and  
by our assumption, the ideal $\big( p(a) \big) \subseteq V$. 
Furthermore, since $a^k=a^kp(a)$ 
(by Proposition \ref{Propo4p(a)}, $iii)$), 
we have $(a^k) \subseteq \big( p(a) \big) \subseteq V$. 
Hence, for all $m\ge k$, we have 
$ba^mc =ba^k(a^{m-k}c)\in (a^k) \subseteq V$. 
\epfv

One immediate consequence of Theorem \ref{CharByIdem} 
is the following corollary which provides a family of 
special $\vt$-Mathieu subspaces.

\begin{corol}\label{CharByIdem-cor1}
Let $V\in \gca$ such that $V$ does 
not contain any nonzero idempotent. 
Then $V$ is a $\vt$-Mathieu subspace 
of $\cA$. 
\end{corol}

\begin{rmk} \label{Def-MS3}
When the algebra $\cA$ is algebraic over $K$, 
by Lemma \ref{3Cases} every $K$-subspace 
$V$ of $\cA$ lies in $\gca$. Then 
Theorem \ref{CharByIdem} gives another 
equivalent formulation for \vmss of 
algebraic $K$-algebras, which is more similar to  
the definition of $\vartheta$-ideals 
than the one given in Definitions \ref{Def-MS}, \ref{Def-MS4} or 
in Proposition \ref{MS-Def2}.
For example, when $\cA$ is algebraic over $K$, 
a  $K$-subspace $M\subseteq \cA$ is a left 
$($resp., right$)$ Mathieu subspace of $\cA$ 
iff for any idempotent $a\in M$ and any $b\in \cA$, 
we have $ba \in M$ $($resp., $ab \in M$$)$. 
\end{rmk}

Next, for any $K$-subspace $V$ of $\cA$ and 
$\vt\neq${\it ``pre-two-sided\hspace{.1mm}"}, we let $I_{\vartheta, V}$ 
denote the $\vartheta$-ideal of $\cA$ which 
is maximum among all the $\vartheta$-ideals 
of $\cA$ contained in $V$. Note that by Zorn's lemma, 
it is easy to see that $I_{\vartheta, V}$ always exists 
and is unique. Actually, $I_{\vartheta, V}$ is the same 
as the sum of all the $\vartheta$-ideals 
of $\cA$ contained in $V$. For example, when $V$ itself is a 
$\vartheta$-ideal of $\cA$, we have $I_{\vartheta, V}=V$. 
In particular, $I_{\vartheta, \cA}=\cA$.

Furthermore, for the case $\vt=$``{\it pre-two-sided\hspace{.1mm}"}, we set 
\begin{align*}
I_{\vartheta, V}\!:=I_{\mbox{\it left},\, V}+ I_{\mbox{\it right},\, V}.
\end{align*}

In other words, $I_{\vartheta, V}$ with  
$\vt=$``{\it pre-two-sided\hspace{.1mm}"} 
is the sum of the maximum left ideal contained in $V$ and 
the maximum right ideal contained in $V$. 
Note that when $\cA$ is not commutative, 
$I_{\vartheta, V}$ in this case is {\it not} necessarily a  
two-sided or one-sided ideal of $\cA$. 

\begin{propo}\label{CharByIdem-cor2}
Let $V\in \gca$ such that $I_{\vartheta, V}=0$. 
Then $V$ is a \vms of $\cA$ iff 
$V$ does not contain any nonzero 
idempotent. 

Consequently, for any $\vt\neq${\it ``pre-two-sided\hspace{.1mm}"} and 
any algebraic $K$-algebra $\cA$ that has no non-trivial  
$\vartheta$-ideals,  
we have that a non-trivial  $K$-subspace $M$ of $\cA$ is 
a \vms of $\cA$ iff $M$ does not contain any nonzero 
idempotent of $\cA$. 
\end{propo} 

\pf The $(\Leftarrow)$ part follows from Corollary 
\ref{CharByIdem-cor1}. To show the 
$(\Rightarrow)$ part, assume that there exists a 
nonzero idempotent $e\in V$. Then 
by Theorem \ref{CharByIdem}, we have 
$(e)_\vartheta \subseteq V$, whence  
$(e)_\vartheta \subseteq I_{\vartheta, V}$. 
Since $0\ne e \in I_{\vartheta, V}$, we have  
$I_{\vartheta, V}\ne 0$, 
which is a contradiction.
\epfv

\begin{corol}
Let $V$ be a  $K$-subspace of $\cA$ and 
$I_V=I_{\vartheta, V}$ with $\vartheta=${\it ``two-sided\hspace{.1mm}"}.  
Assume that $V\in \gca$ or 
$V/I_{V}\in \mathcal G(\cA/I_V)$. 
Then $V$ is a Mathieu subspace of $\cA$ iff $V /I_{V}$ 
does not contain any nonzero idempotent 
of the quotient $K$-algebra $\cA/I_{V}$.
\end{corol}
\pf First, 
%since $I_V\subseteq V$, 
it is easy to see that 
$V\in \gca$ implies $V/I_V \in \mathcal G(\cA/I_V)$.
So we may assume the latter.  
Second, since $I_V$ is maximum among all the   
ideals of $\cA$ that are contained in $V$, 
the quotient $V /I_V$ does not contain 
any nonzero ideal of the quotient algebra $\cA/I_V$, 
whence $I_{V/I_V}=0$. 

Now, applying Proposition \ref{CharByIdem-cor2} 
to the $K$-algebra $\cA/I_V$ and 
its $K$-subspace $V /I_V$, we see that 
$V /I_V$ is a Mathieu subspace of 
$\cA/I_V$ iff $V /I_V$ does not contain 
any nonzero idempotent of $\cA/I_V$. 
On the other hand, by Proposition \ref{quotient-I} 
we also have that 
$V$ is a Mathieu subspace of $\cA$ iff 
$V /I_V$ is a Mathieu subspace of 
$\cA/I_V$. Combining these 
two equivalences the corollary follows.
\epfv

Next we derive some consequences of 
Corollary \ref{CharByIdem-cor1} on 
finite dimensional \vmss of 
$K$-algebras. 

\begin{propo}\label{FinDim}
Assume that $\cA$ is purely transcendental 
over $K$, i.e., the only algebraic elements of 
$\cA$ are the elements in $K\subseteq \cA$. 
Then every finite dimensional $K$-subspace 
$V$ of $\cA$ such that $1\not \in V$ 
is a $\vt$-Mathieu subspace of $\cA$.  
\end{propo}

\pf Since $\cA$ is purely transcendental 
over $K$ and all idempotents of $\cA$ 
are algebraic over $K$, we see that 
all idempotents of $\cA$ must 
lie inside $K\subseteq \cA$. 

But, on the other hand, all idempotents 
of $K$ are the solutions of the equation 
$t^2-t=0$ in $K$, which are $0, 1\in K$.  
Therefore, all idempotents of $\cA$  
are trivial. Furthermore, 
since $1\not \in V$, we see that $V$ does not 
contain any nonzero idempotent of $\cA$. 
Then the proposition follows 
immediately from Lemma \ref{3Cases} and 
Corollary \ref{CharByIdem-cor1}.  
\epfv 

The following characterization of one-dimensional 
\vmss of associative $K$-algebras will play important 
roles in the later Sections \ref{S5}-\ref{S7}. 

\begin{propo}\label{1D-MS}
Let $\cA$ be an associative $K$-algebra and 
$0\ne a\in \cA$. 
Then the one-dimensional $K$-subspace $Ka$ 
is a \vms of $\cA$ iff one of the following 
two statements holds:

$1)$ $Ka$ is a $\vartheta$-ideal of $\cA$, 
or equivalently, $Ka=(a)_\vartheta$.

$2)$ $a$ is not a quasi-idempotent of $\cA$. 
\end{propo}

It is an easy exercise to check that when 
$\vt=${\it ``pre-two-sided\hspace{.1mm}"}, 
the equivalence in $1)$ above indeed holds, i.e., 
$Ka$ is a {\it pre-two-sided} ideal of $\cA$, 
which by definition means a ({\it two-sided}) ideal, 
iff $Ka=(a)_\vartheta=aA+Aa$. \\

\underline{\it Proof of Proposition \ref{1D-MS}:} \, 
$(\Rightarrow)$ Assume that $Ka$ is a \vms of $\cA$  
but statement $2)$ fails, i.e., $a$ is a nonzero  
quasi-idempotent of $\cA$. 
%
%Then we have 
%$a^2 =r a$ for some $r\in K^\times$. 
%
%Let $e\!:=r^{-1}a$. Then it is easy to check that 
%$e$ is an idempotent and lies in $Ka$. 
Then by Lemma \ref{CyclicLemma} $iii)$, we have 
$(a)_\vartheta \subseteq Ka$. 
Since $(a)_\vartheta \supseteq Ka$, 
we have $(a)_\vartheta = Ka$, i.e., 
statement $1)$ holds.
 
$(\Leftarrow)$ If statement $1)$ holds, 
then $Ka$ is a $\vartheta$-ideal of $\cA$  
and hence, also a \vms of $\cA$. 
Assume that statement $2)$ holds. 
Then for any $r\in K^\times$, 
$b\!:=ra$ cannot be an idempotent, 
otherwise $a=r^{-1}b$ would be a quasi-idempotent too. 
Hence, $Ka$ does not contain any nonzero idempotent 
of $\cA$. Then by Lemma \ref{3Cases} and Corollary \ref{CharByIdem-cor1}, 
$Ka$ is a \vms of $\cA$. 
\epfv

\subsection{Radicals of $\vartheta$-Mathieu Subspaces $M\in \ovca$ 
in Terms of Radicals of $I_{\vartheta, M}$}\label{S4.2}

Throughout this subsection, for each \vms $M$ of $\cA$, 
for convenience we denote by $I_M$ 
the notation $I_{\vartheta, M}$ introduced in 
the previous subsection. In particular, when 
$\vt\neq${\it ``pre-two-sided\hspace{.1mm}"},   
$I_M$ denotes the unique $\vartheta$-ideal 
of $\cA$ which is maximum among all the 
$\vartheta$-ideals of $\cA$ contained in $M$.

\begin{lemma}\label{sqM=sqI}
Let $\cA$ be a  $K$-algebra and 
$M$ a \vms of $\cA$. Then $\arad(M)=\arad(I_M)$.

In particular, if $M\in \ovca$, we have 
$\sqrt M=\sqrt{I_M}$.
\end{lemma}

\pf Note first that since $M\supseteq I_M$, we have 
$\arad(M) \supseteq \arad(I_M)$. 

To show $\arad(M) \subseteq \arad(I_M)$, 
let $a\in \arad(M)$. 
Since $a$ is algebraic over 
$K$, it follows from Theorem \ref{Char4AlgElt} 
that $(a^N)_\vartheta\subseteq M$ for some 
$N\ge 0$. Hence, we also have 
$(a^N)_\vartheta\subseteq I_M$. 
Consequently,  $a^m\in I_M$ 
for all $m\ge N$, whence $a\in \arad(I_M)$.
\epfv

%\begin{lemma}\label{sqM=sqI}
%For any $M\in \ovca$, we have 
%$\sqrt M=\rad({I_{\vartheta, M}})$.
%\end{lemma}
%
%\pf Note first that since $M\supseteq I_M$, we have 
%$\sqrt M\supseteq \sqrt {I_M}$. 
%
%To show $\sqrt M\subseteq \sqrt {I_M}$, 
%pick up any $a\in \sqrt M$. 
%Since $a$ is algebraic over 
%$K$, it follows from Theorem \ref{Char4AlgElt} 
%that $(a^N)_\vartheta\subset M$ for some 
%$N\ge 1$. Hence, we also have 
%$(a^N)_\vartheta\subset I_M$. 
%Consequently, we have $a^m\in I_M$ 
%for all $m\ge N$, whence $a\in \sqrt{I_M}$.
%\epfv

\begin{theo}\label{sqM=sqI-Thm}
Let $M\in \ovca$ and $V$ a  $K$-subspace of $M$ 
such that $I_M\subseteq V$. Then $V$ is a \vms of 
$\cA$ and $\sqrt V =\sqrt{I_M}$. 
\end{theo}

\pf Note first that by Lemma \ref{sqM=sqI}, 
it suffices to show that $V$ is a \vms of $\cA$, for  
we obviously have $V\in \gca$ and $I_V=I_M$. 

Let $e$ be a nonzero idempotent in $V$. 
Hence, $e\in M$ since $V\subseteq M$. 
Then by Theorem \ref{CharByIdem}, we have 
$(e)_\vartheta \subseteq M$, whence 
$(e)_\vartheta \subseteq I_M \subseteq V$. 
Then by Theorem \ref{CharByIdem} again, 
$V$ is a \vms of $\cA$. 
\epfv

\begin{corol}\label{SimpleAlg-2}
Let $\cA$ be a simple and algebraic $K$-algebra 
and $M$ a proper Mathieu subspace of $\cA$. Then 
$\sqrt M = \nil(\cA)$ and 
all $K$-subspaces $V\subseteq M$ 
are also Mathieu subspaces of $\cA$. 
\end{corol}

\pf Since $\cA$ is simple, we have $I_M=0$. Since $\cA$ 
is algebraic over $K$, by Lemmas \ref{3Cases} and \ref{sqM=sqI} 
we have $\sqrt M=\sqrt 0=\nil(\cA)$. Then  
the corollary follows from Theorem \ref{sqM=sqI-Thm} 
or Lemma \ref{radical-Lemma}.
\epfv

When the $K$-algebra $\cA$ is commutative, 
we have the following characterization 
for the Mathieu subspaces with 
algebraic radicals.

\begin{theo}\label{Comm-MSs}
Let $\cA$ be a commutative $K$-algebra and 
$V\in \gca$. Then $V$ is a Mathieu subspace of $\cA$ 
iff $\sqrt V$ is an ideal of $\cA$.
\end{theo}

\pf It is well-known that the radicals of 
ideals of commutative algebras are  
(radical) ideals. Then the $(\Rightarrow)$ 
part follows immediately 
from Lemma \ref{sqM=sqI}, 
for $\sqrt V=\sqrt{I_V}$ and 
$I_V$ is an ideal of $\cA$. 

To show the $(\Leftarrow)$ part, 
by Theorem \ref{CharByIdem} it suffices 
to show that for each idempotent $e\in V$, 
we have $(e) \subseteq V$. Equivalently, 
it suffices to show that
the $K$-subspace $V_e \!:=\{a\in \cA\,|\, ea \in V\}$ 
is equals to $\cA$ itself. 

Note first that by Lemma \ref{CyclicLemma}, $ii)$ 
we have $e\in \sqrt V$. 
Since $\sqrt V$ by our assumption is an ideal of $\cA$,  
we have $eb\in \sqrt V$ for all $b\in \cA$. 
Then for all $m\gg 0$, we have $eb^m=(eb)^m\in V$ 
or equivalently, $b^m\in V_e$. Hence,  
$b\in \sqrt{V_e}$ for all $b\in \cA$, 
whence $\sqrt{V_e}=\cA$. Applying  
Lemma \ref{CuteLemma} to the 
$K$-subspace $V_e$, we get $V_e=\cA$.
\epfv

One by-product of Theorem \ref{Comm-MSs} is the 
following corollary which does not seem obvious.

\begin{corol}\label{Comm-MSs-Corol}
Let $\cA$ be a commutative $K$-algebra and 
$V\in \gca$. Then $\sqrt V$ is a radical ideal of $\cA$ 
if $($and only if$)$ $\sqrt V$ is an ideal of $\cA$.
\end{corol}
\pf Assume that $\sqrt V$ is an ideal of $\cA$. Then  
by Theorem \ref{Comm-MSs}, $V$ is a Mathieu subspace of $\cA$, and 
by Lemma \ref{General-RadLemma}, $\sqrt V$ is a 
radical ideal of $\cA$.
\epfv

Next, we conclude this subsection with the following two remarks.  

%\vspace{2mm}

First, as we can see from the example below, without the algebraic condition on $\sqrt V$ Theorem \ref{Comm-MSs} does not always hold. 

%One remark on Theorem \ref{Comm-MSs} is that it 
%does not always hold without the algebraic 
%condition on $\sqrt V$. 

\begin{exam}
Let $\cA$ be the Laurent polynomial algebra $\bC[t^{-1}, t]$ 
in one variable $t$ over $\bC$ and $V$ the subspace of 
all Laurent polynomials in $\cA$ without constant terms.
Then by the Duistermaat-van der Kallen theorem \cite{DK}, $V$ is a Mathieu subspace of $\cA$ and $\sqrt V=t\bC[t]\bigcup t^{-1}\bC[t^{-1}]$, which is not even a  $\bC$-subspace and hence, 
not an ideal of $\cA$. 
\end{exam}

Second, even though the univariate polynomial algebra $K[t]$ 
is purely transcendental over $K$, by using 
Theorems \ref{sqM=sqI-Thm} and \ref{Comm-MSs},  
it has been shown recently in \cite{EZ} that 
the following theorem actually also holds. 
 
\begin{theo} $($\cite{EZ}$)$\, 
Let $V$ be a $K$-subspace of the univariate 
polynomial algebra $K[t]$. 
Then $V$ is a Mathieu subspace of $K[t]$ 
iff $\sqrt V=\sqrt{I_V}$. 
\end{theo}

But, it  has also been shown in \cite{EZ} that the theorem above fails 
for multi-variable polynomial algebras.

%{\bf Question:} But it is 
%interesting to see if the $(\Leftarrow)$ part 
%of the theorem still holds without the 
%algebraic condition on $\sqrt V$.

%\begin{corol}
%Let $\cA$ be a commutative $K$-algebra and 
%$V$ any $K$-subspace of $\cA$ such that 
%$\sqrt V$ is algebraic. 
%Assume that $\sqrt M$ is not a radical 
%ideal of $\cA$. Then $M$ is not a Mathieu subspace of $\cA$.
%\end{corol}

\subsection{Unions and Intersections of Mathieu Subspaces with Algebraic Radicals} \label{S4.3}

First, let's prove the following proposition on the 
intersections of $\vartheta$-Mathieu subspaces.

\begin{propo}\label{Interscn}
Let $M_i$ $(i\in I)$ be a family of \vmss of 
a  $K$-algebra $\cA$. Assume that $M_i\in \gca$ for some 
$i\in I$, or the intersection 
$\bigcap_{i\in I} M_i\in \gca$. 
Then $\bigcap_{i\in I} M_i$ is also a \vms of $\cA$.
\end{propo}

\pf Note first that the condition that 
$M_i\in \gca$ for some $i\in I$ obviously implies 
the condition $\bigcap_{i\in I} M_i\in \gca$. 
So we may assume the latter. 

Let $e$ be any idempotent in $\bigcap_{k\ge 1} M_k$. 
Then for any $i\in I$, we have $e\in M_i$ and by 
Lemma \ref{CyclicLemma} $iii)$,   
$(e)_\vartheta \subseteq M_i$, whence 
$(e)_\vartheta \subseteq \bigcap_{i\in I} M_i$.  
Applying Theorem \ref{CharByIdem} to  
$\bigcap_{i\in I} M_i$, 
the proposition follows. 
\epfv

It is worthy to point out that it is easy to check 
(or see Proposition $4.9$ in \cite{GIC}) that 
in general the intersection of any finitely 
many \vmss is always a $\vt$-Mathieu subspace.   
However, when $|I|=\infty$, Proposition \ref{Interscn} 
without the algebraic conditions does not always hold.

\begin{exam}\label{Counter-Exam1}
Let $\cA$ be the polynomial algebra $K[t]$ 
in one variable $t$ over $K$ and 
$M_i$ $(i\ge 0)$ the $K$-subspace 
of $\cA$ spanned by the monomials $t^k$ 
with $k\ge 1$ but $k\ne 2j+1$ for 
all $0\le j\le i$. Then 
it is easy to check that for each $i\ge 0$, 
$\sqrt{M_i}=tK[t]$ and $M_i$ is a 
Mathieu subspace of $K[t]$.

On the other hand, we also have 
$M\!:=\bigcap_{i\ge 0}M_i=t^2K[t^2]$. 
Note that $t^2\in \sqrt M$. But,  
for each $m\ge 1$, $t(t^2)^m=t^{2m+1}\not \in M$.
Hence, the intersection $M$ of $M_i$ $(i\ge 0)$ 
is not a Mathieu subspace of $\cA$.
\end{exam}

Next we consider unions of ascending sequences 
of \vmss under certain conditions. 
  
\begin{propo}\label{Union}
Let $M_i$ $(i\ge 1)$ be a sequence of non-trivial  \vmss of $\cA$ 
such that $M_i\subseteq M_{i+1}$ for all $i\ge 1$.
Assume that $\bigcup_{i\ge 1} M_i \in \gca$.  
Then $\bigcup_{i\ge 1} M_i$ is also a non-trivial  \vms of $\cA$.
\end{propo}

\pf First, since $M_i$ is non-trivial  for each $i\ge 1$, 
by Lemma \ref{OneLemma} we have 
$1\not \in M_i$, whence $1\not \in \bigcup_{i\ge 1} M_i$.
So the union $\bigcup_{i\ge 1} M_i$ is also non-trivial. 

Second, let $e$ be an idempotent in $\bigcup_{i\ge 1} M_i$.  
Then $e\in M_k$ for some $k\ge 1$, and by    
Lemma \ref{CyclicLemma} $iii)$, $(e)_\vartheta\subseteq M_k$.
Hence we have $(e)_\vartheta \subseteq \bigcup_{i\ge 1} M_i$.  
Then by Theorem \ref{CharByIdem},    
$\bigcup_{i\ge 1} M_i$ is a \vms of $\cA$. 
\epfv

One remark on Proposition \ref{Union} is that 
without the algebraic condition on the 
radical of the union, the proposition  
does not necessarily hold. 

\begin{exam}
Let $\cA=K[t]$ as in Example \ref{Counter-Exam1} 
and $V_i$ $(i\ge 1)$ the $K$-subspace of $\cA$ 
spanned by the monomials $t^{2j}$ $(1\le j\le i)$.
Note that for any $i\ge 1$, we have 
$1\not \in V_i$ and $\dim_K V_i<\infty$.  
Then it follows from Proposition \ref{FinDim} that  
for any $i\ge 1$, $V_i$ is a 
Mathieu subspace of $\cA$.

But, on the other hand, we have   
$\bigcup_{i\ge 1} V_i=t^2K[t^2]$, 
which as shown in Example \ref{Counter-Exam1},  
is not a Mathieu subspace of $\cA$. 
\end{exam}

Next, we use Zorn's lemma and 
Propositions \ref{Interscn} and \ref{Union} to derive 
existences of certain {\it maximal} (resp., {\it minimal}) 
non-trivial $\vartheta$-Mathieu subspaces for 
algebraic $K$-algebras $\cA$. 

%But, before we proceed, 
%let's fix the following 
%general definition.
%
%\begin{defi}
%Let $\cA$ be a $R$-algebra and $V$ a $R$-subspace of $\cA$. 
%We say $V$ is a {\it maximal} $($resp., {\it minimal}$)$ 
%{\it non-trivial} \vms of $\cA$ if $V$ is a maximal 
%$($resp., {\it minimal}$)$ element in the collection 
%of all non-trivial \vmss of $\cA$. 
%\end{defi}
%
%Of course, the maximal and minimal 
%{\it non-trivial} \vms may not exist, 
%and in general very little is known about 
%these extremal non-trivial $\vartheta$-Mathieu 
%subspaces. But for algebraic algebras over  
%fields, as we will see below the situation 
%becomes much better.  

First, note that if $\cA$ is 
algebraic over $K$, then by Lemma \ref{3Cases} 
the algebraic conditions in Propositions \ref{Interscn} 
and \ref{Union} are automatically 
satisfied. With this observation and by Zorn's lemma, 
we immediately have the following 
proposition. 

\begin{propo}\label{MaxMinMSs-Propo}
Let $\cA$ be an algebraic $K$-algebra 
and $V$ a  $K$-subspace of $\cA$. 
Then the following statements hold.  

$i)$ There exists at least one \vms of $\cA$
which is maximal among all the 
\vmss of $\cA$ contained in $V$.

$ii)$ There exists a unique \vms $M$ of $\cA$ 
which is minimum among all the  
\vmss $W$ of $\cA$ with 
$V\subseteq W$. Actually, $M$ 
is given by the intersection 
of all \vmss that contain $V$.  

$iii)$ Any non-empty collection 
of proper \vmss $M$ of $\cA$ with 
$V\subseteq M$ has at least one maximal element  
and a $($unique$)$ minimum element. 
\end{propo}

\begin{theo}\label{MaxMSs-Thm}
Assume that $\cA$ is algebraic over $K$ but   
$\cA\ne K$.  Then for any proper \vms $M$ of $\cA$, 
there exists a maximal non-trivial \vms of $\cA$ 
which contains $M$.

In particular, $($by taking $M=0$$)$, 
$\cA$ has at least one maximal non-trivial   
$\vartheta$-Mathieu subspace.
\end{theo}

\pf Let $\cF$ be the collection of the non-trivial  \vmss 
$J$ of $\cA$ such that $M\subseteq J$. 
If $M\ne 0$, then $M\in \cF$. 
If $M=0$, then by Lemma \ref{FieldCase} in the later 
Section \ref{S6}, $\cA$ has at least one non-trivial  
$\vartheta$-Mathieu subspace $J$, 
which obviously lies in $\cF$. 
Therefore, in any case  
$\cF\ne \emptyset$. Then the theorem follows 
directly from Proposition \ref{MaxMinMSs-Propo} 
$iii)$.
\epfv
  
\begin{corol}\label{MaxMSs-Corol}
Let $V$ be a $K$-subspace of an algebraic $K$-algebra $\cA$ 
such that the $\vartheta$-ideal generated by elements of 
$V$ is non-trivial . Then there exists a maximal non-trivial  
\vms $M$ of $\cA$ such that $V\subseteq M$.
\end{corol}
\pf Since any $\vartheta$-ideal is a 
$\vartheta$-Mathieu subspace, the corollary follows 
immediately from Theorem \ref{MaxMSs-Thm} by taking 
$M$ to be the $\vartheta$-ideal generated 
by elements of $V$.
\epfv

\renewcommand{\theequation}{\thesection.\arabic{equation}}
\renewcommand{\therema}{\thesection.\arabic{rema}}
\setcounter{equation}{0}
\setcounter{rema}{0}

\section{\bf Co-dimension One Mathieu Subspaces and the Minimal Non-trivial 
Mathieu Subspaces of Matrix Algebras over Fields}\label{S5}

Let $K$ be an arbitrary field and $n\ge 1$. 
In this section we classify the co-dimension 
one \vmss and the minimal non-trivial \vmss  
for the matrix algebra $M_n(K)$ of 
$n\times n$ matrices with entries in $K$.

First, let's fix the following notations that will be used throughout this section.

We denote by $I_n$ the {\it identity} matrix in $M_n(K)$. For each $X\in M_n(K)$, we denote by $\Tr X$ 
the {\it trace} of the matrix $X$ and set
\begin{align}\label{Def-NX} 
H_X\!:=\{A\in M_n(K)\,|\, \Tr (AX)=0 \}.
\end{align}
When $X=I_n$, $H_{I_n}$ will also be denoted by $H$, 
i.e.,  
\begin{align}\label{Def-H} 
H\!:=\{A\in M_n(K)\,|\, \Tr A=0\}.
\end{align}

For any $X, Y\in M_n(K)$, we denote by $X\sim Y$ 
if $X=sY$ for some $s\in K^\times$. Note that 
by Lemma \ref{Co-D1-Lemma1} below,  
we have 
\begin{align}
H_X=H_Y \, \Leftrightarrow \, X\sim Y. \label{NX=NY} 
\end{align}
In particular, we have  
\begin{align}
H_X=H \, \Leftrightarrow \, X \sim I_n. \label{NX=H} 
\end{align}

With the notations fixed above, the first main result of this section can be stated as follows.

\begin{theo}\label{Co-D1}
Let $K$ be a field and $n\ge 1$. Then the following two statements hold: 

$i)$ if $char.\,K=0$ or $char.\,K=p>n$, then 
$H$ is the only co-dimension one 
\vms of $M_n(K)$; 

$ii)$ if $char.\,K=p>0$ and $p\le n$, 
then $M_n(K)$ has no co-dimension one 
$\vartheta$-Mathieu subspaces.
\end{theo}

In order to prove the theorem, 
we first need to prove 
the following two lemmas.

\begin{lemma}\label{Co-D1-Lemma1}
For every co-dimension one $K$-subspace 
$V$ of $M_n(K)$, there exists $0\ne X\in M_n(K)$ 
such that $V=H_X$. Furthermore, 
$X$ is unique up to nonzero 
scalar multiplications. 
\end{lemma}
\pf First, let's consider 
the following $K$-bilinear form of 
$M_n(K)$: 
\begin{align}
(\cdot, \cdot): M_n(K)\times M_n(K) &\to \quad  K  
\label{T-Pairing} \\
(A,\quad B)\quad \quad &\to \,  \Tr(AB). \nno 
\end{align}

It is well-known and also easy to check that 
the bilinear form above is non-singular. 
Hence, it induces a $K$-linear isomorphism 
\begin{align}
\phi: M_n(K)& 
\overset{\sim}{\rightarrow} \hom_K (M_n(K), K) 
\label{K-Iso} \\
B\quad  &\rightarrow \quad\quad\quad\quad \phi_B,  
\nno
\end{align}
where $\phi_B: M_n(K)\to K$ is the linear functional of 
$M_n(K)$ defined by setting for all $A\in M_n(K)$,    
\begin{align}\label{Def-phi}
\phi_B(A)\!:=\Tr(AB).
\end{align}

Note that any co-dimension one $K$-subspace 
of $M_n(K)$ is the kernel of a nonzero linear 
functional of $M_n(K)$, which is unique up 
to nonzero scalar multiplications. 
Then by the $K$-linear isomorphism 
in Eq.\,(\ref{K-Iso}), we see that 
for the co-dimension one subspace 
$V$ in the lemma, there exists 
$0\ne X\in M_n(K)$, which is unique up 
to nonzero scalar multiplications, 
such that $V=\Ker \, \phi_X$. Furthermore, 
by Eqs.\,(\ref{Def-phi}) and (\ref{Def-NX}), 
we also have $\Ker \, \phi_X=H_X$,  
whence the lemma follows. 
\epfv

\begin{lemma}\label{MnCase-Lemma}
Let $n\ge 2$ and $0\ne X \in M_n(K)$ such that 
$X\not\sim I_n$. Then there exist 
non-trivial idempotents  
$A, B \in M_n(K)$ such that 
\begin{align} 
AX & \ne 0; \label{MnCase-Lemma-e1}\\
XB & \ne 0;\label{MnCase-Lemma-e2} \\ 
\Tr (AX)&=\Tr(XB)=0.\label{MnCase-Lemma-e3}
\end{align}
\end{lemma}

\pf First, it is easy to check that 
the existence of the idempotent 
$B$ for $X$ follows from that of 
the idempotent $A$ for $X^\tau$ 
by letting $B=A^\tau$, 
where $X^\tau$ and $A^\tau$ are  
the transposes of $X$ and $A$, 
respectively. So it suffices to show 
the existence of the non-trivial 
idempotent $A$.

We first consider the case $n=2$. 
Write 
$X=\begin{pmatrix} a& b\\
c&d
\end{pmatrix}
$
for some $a, b, c, d\in K$. 
We divide the proof 
into the following three 
different cases. 

\underline{\it Case 1}\,: If $b\ne 0$, let 
$A\!:=\begin{pmatrix} 1 & 0\\
-ab^{-1} & 0
\end{pmatrix}$.
Then we have 
\begin{align}
AX&=\begin{pmatrix} a & b\\
-a^2b^{-1} & -a
\end{pmatrix}
\ne 0
\label{Case1-e1} \\
\Tr(AX)& =0\quad  \mbox{ and } \quad A^2=A. 
\end{align}

\underline{\it Case 2}\,: If $b=0$ but $c\ne 0$, let 
$A\!:=\begin{pmatrix} 0& -c^{-1}d \\
0 & 1
\end{pmatrix}$.
Then we have 
\begin{align}
AX&=\begin{pmatrix} -d & -c^{-1}d^2 \\
c & d
\end{pmatrix} 
\ne 0 \\
\Tr(AX)&=0 \quad \mbox{ and } \quad A^2=A. \label{Case2-e2}
\end{align}

\underline{\it Case 3}\,: If $b=c=0$, 
then $a\ne d$ since by our assumption 
$X\not \sim I_2$. In particular, 
$a$ and $d$ cannot 
be both zero.

Let $A\!:=\frac{1}{d-a}\begin{pmatrix} d& d \\
-a & -a
\end{pmatrix}$.
Then we have 
\begin{align}
AX&=\frac{1}{d-a}\begin{pmatrix} ad & d^2 \\
-a^2 & -ad 
\end{pmatrix} 
\ne 0
\label{Case3-e1} \\
&\Tr(AX)=0 \quad \mbox{ and } \quad A^2=A. \label{Case3-e2}
\end{align}

It is straightforward to check that all 
the equations (\ref{Case1-e1})--(\ref{Case3-e2})  
do hold and that the idempotent $A$ in each case is non-trivial. So we omit the details here.

Next, we consider the case $n\ge 3$. 
Since $X\not \sim I_n$, it is easy to 
see that there exist $1\le m<k\le n$ 
such that the $2\times 2$-minor of $X$ on 
$m^{\rm th}$, $k^{\rm th}$ 
rows and $m^{\rm th}$, $k^{\rm th}$ 
columns is not a multiple 
of $I_2$. 

Since in general idempotents and also traces of 
matrices are preserved by conjugations,  
by applying some conjugations by permutation 
matrices to $X$ if it is necessary, 
we may further assume $m=1$ and $k=2$. 
We denote by $X'$ this $2\times 2$ 
minor of $X$. 

By the lemma for the case $n=2$, 
there exists a non-trivial idempotent  
$A'\in M_2(K)$ such that 
$A'X'\ne 0$ and $\Tr (A'X')=0$.
Let 
$A\!:=\begin{pmatrix} A' & 0 \\
0 & 0 
\end{pmatrix}
\in M_n(K)$. 
Then it is easy to check  
that $A$ is a non-trivial idempotent 
of $M_n(K)$ which satisfies 
$AX\ne 0$ and $\Tr (AX)=0$.
Hence, the lemma also holds 
for the case $n\ge 3$.
\epfv

Note that one bi-product of Lemmas \ref{Co-D1-Lemma1} 
and \ref{MnCase-Lemma} above is the following corollary. 

\begin{corol}\label{MnCase-Corol}
Let $n\ge 2$ and $V$ be 
a co-dimension one $K$-subspace of $M_n(K)$  
such that $V\ne H$. Then $V$ contains at least 
one non-trivial idempotent of $M_n(K)$. 
\end{corol}
\pf By Lemma \ref{Co-D1-Lemma1}, we know that 
$V=H_X$ for some $0\ne X\in M_n(K)$. Since $V\ne H$, we have $X\not \sim I_n$. Then by  
Lemma \ref{MnCase-Lemma}, there exists a  
non-trivial idempotent  $A$ of 
$M_n(K)$, which satisfis 
Eq.\,(\ref{MnCase-Lemma-e3}). 
Hence by Eq.\,(\ref{Def-NX}), 
we have $A, B\in H_X=V$. 
\epfv 

Now we can prove the first main result of this section as follows.\\

\underline{\it Proof of Theorem \ref{Co-D1}:} \, 
Note first that if $n=1$, then the theorem obviously holds. 
So we assume $n\ge 2$. 

Let $V$ be a co-dimension one $K$-subspace 
of $M_n(K)$ such that $V\ne H$. 
Then by Lemma \ref{Co-D1-Lemma1}, 
$V=H_X$ for some $0\ne X\in M_n(X)$. 
Note that by Eq.\,(\ref{NX=H}), 
$X\not \sim I_n$ since $V\ne H$. 
Next, we show that $V$ cannot 
be a left or right Mathieu 
subspace of $M_n(X)$.

Assume that $V$ is a left  Mathieu subspace 
of $M_n(K)$. Let $A$ be the non-trivial 
idempotent as in Lemma \ref{MnCase-Lemma}. 
Then by Eqs.\,(\ref{Def-NX}) and 
(\ref{MnCase-Lemma-e3}), we have 
$A\in H_X=V$. Furthermore, 
by Lemma \ref{CyclicLemma} $iii)$, 
we have $CA\in V=H_X$ for all $C\in M_n(K)$. 
More precisely, 
we have
\begin{align}\label{MarkEq}
0=\Tr \big((CA)X\big)=\Tr\big(C(AX)\big)
\end{align}
for all $C\in M_n(K)$.

Since the $K$-bilinear form in 
Eq.\,(\ref{T-Pairing}) is 
non-singular, we have $AX=0$. 
But this contradicts  
Eq.\,(\ref{MnCase-Lemma-e1}) 
in Lemma \ref{MnCase-Lemma}. 
Therefore, $V$ cannot be a 
left Mathieu subspace of 
$M_n(K)$. 

Assume that $V$ is a right  Mathieu subspace 
of $M_n(K)$. Let $B$ be the non-trivial 
idempotent as in Lemma \ref{MnCase-Lemma}. 
Then by Eqs.\,(\ref{Def-NX}) and 
(\ref{MnCase-Lemma-e3})  we have 
$B\in H_X=V$, and by Lemma \ref{CyclicLemma} $iii)$, 
$BC\in V=H_X$ for all $C\in M_n(K)$. 
More precisely, we have 
\begin{align*}
0=\Tr \big((BC)X\big)=\Tr \big(B(CX)\big)
=\Tr\big((CX)B\big)=\Tr\big(C(XB) \big)
\end{align*}
for all $C\in M_n(K)$.

Then by the non-singularity of the 
$K$-bilinear form in 
Eq.\,(\ref{T-Pairing}) again, 
we have $XB=0$, which contradicts  
Eq.\,(\ref{MnCase-Lemma-e2}) 
in Lemma \ref{MnCase-Lemma}. 
Therefore, $V$ cannot be a 
right Mathieu subspace of 
$M_n(K)$ either. 

Therefore, for any specification of $\vt$, the only possible 
co-dimension one \vms of $M_n(K)$ is the $K$-subspace $H$ of the 
trace-zero matrices in $M_n(K)$, which we will consider next.

Assume first $char.\,K=p \le n$. 
Let $e_p\!:=\begin{pmatrix} 
I_p & 0 \\
0 & 0
\end{pmatrix}
\in M_n(K).$ 
Note that $e_p$ is a nonzero idempotent lying in $H$, and 
$(e_p)_\vartheta$ clearly contains the subalgebra 
$\begin{pmatrix} 
M_p(K) & 0 \\
0 & 0
\end{pmatrix}\subseteq M_n(K)$, 
which certainly cannot be entirely contained in $H$. 
Therefore, we have $(e_p)_\vartheta\not\subseteq H$. 
Then by Lemma \ref{CyclicLemma}, $iii)$ or Theorem \ref{CharByIdem}, 
$H$ in this case cannot be a $\vartheta$-Mathieu subspace 
of $M_n(K)$, whence the statement $ii)$ 
of the theorem follows.
 
%Then by the similar 
%arguments as above for the case $V=H_X$ 
%$(X\not \sim I_n)$ with $H_X$ 
%replaced by $H=H_{I_n}$ 
%(i.e., replacing $X$ by $I_n$)  
%and the idempotents $A$ and $B$  
%by the idempotent $D_p$, 
%we see that $H$ in this case 
%is not a left or right 
%Mathieu subspace of $M_n(K)$ and hence, 
%not a (two-sided) Mathieu subspace either. 
%Therefore, the statement $ii)$ 
%in the theorem holds. 

Now, assume $char.\,K=0$ or $char.\,K=p>n$. 
Then it is well-known in linear algebra that 
for any $A\in M_n(K)$, $A$ is nilpotent iff 
for all $m\ge 1$, $\Tr (A^m)=0$, i.e., $A^m\in H$. 
Hence, we have $\sqrt{H}=\nil(M_n(K))$. Then by 
Lemma \ref{radical-Lemma} the statement $i)$ of 
the theorem also follows. 
\epfv
%
%Next we conclude this section with the following 
%remarks on the minimal and maximal  
%$\vt$-Mathieu subspaces of 
%$M_n(K)$.\\
%
%First, from Proposition \ref{1D-MS}  
%it is easy to see that $M_n(K)$ 
%$(n\ge 2)$ always has non-trivial  
%$\vt$-Mathieu subspaces, for not all matrices 
%of $M_n(K)$ can be quasi-idempotents. 
%Since $M_n(K)$ is of finite dimension, 
%hence it has minimal and maximal non-trivial  
%$\vt$-Mathieu subspaces. 
%For the minimal non-trivial  ones, 
%we have the following classification.  

Next we give a classification for the minimal 
non-trivial \vmss of the matrix algebras 
$M_n(K)$ $(n\ge 2)$. 

\begin{propo}\label{ClsMin}
A $K$-subspace $V \subset M_n(K)$ $(n\ge 2)$ is 
a minimal non-trivial $\vt$-Mathieu subspace of  
$M_n(K)$ iff $V=KA$ for some nonzero 
$A\in M_n(K)$ which is not a 
quasi-idempotent.         
\end{propo}

%\pf First, it is well-known that $M_n(K)$ 
%is a simple $K$-algebra (e.g., see the lemma 
%on p.\,$9$ in \cite{P}). Second, $M_n(K)$ 
%is also algebraic over $K$ since it is of 
%finite dimension. 
%
%Now, let $V$ be a minimal non-trivial   
%Mathieu subspace of $M_n$. Then 
%by Corollary \ref{SimpleAlg-2}, 
%every $K$-subspace of $V$ is also 
%a Mathieu subspace of $M_n(K)$. 
%Since $V$ is nonzero and minimum, 
%we have $\dim_K V=1$. Then by Proposition \ref{1D-MS} 
%and the fact that $M_n(K)$ is a simple $K$-algebra, 
%it is easy to see that the one-dimensional 
%$K$-subspace $V$ cannot be spanned by 
%any idempotent of $M_n(K)$. 
%\epfv

To prove the proposition, we need first to show the following lemma. 

\begin{lemma}\label{ClsMin-L} 
For any $n\ge 2$ and $0\ne A\in M_n(K)$, we have  
%and $J$ the left or right ideal of $M_n(K)$ generated by $A$. 
\begin{enumerate}
  \item[$i)$] $(A)_\vt \ne KA$; 
  \item[$ii)$] $(A)_\vt$ contains at least 
  one element which is not a quasi-idempotent.  
\end{enumerate} 
\end{lemma}

Note that from the well-known fact that $M_n(K)$ 
is a simple $K$-algebra (e.g., see the lemma 
on p.\,$9$ in \cite{P}), it follows immediately 
that the lemma holds for the two-sided case, 
since in this case $(A)_\vt=(A)=M_n(K)$. 
But, for the other cases, we need a different 
argument given below, which actually works 
for all the cases. \\

\underline{\it Proof of Lemma \ref{ClsMin-L}:}\, 
Note first that for any $\vt$, $(A)_\vt$  
contains either the left ideal generated by 
$A$ or the right ideal generated by $A$. 
Therefore, it suffices to show the proposition 
for the two cases: $\vt=${\it ``left"} and 
$\vt=${\it ``right"}. 

We here just give a proof for the former case.   
The latter case follows from 
the former one for the transpose $A^\tau$ of $A$,  
or by applying the similar arguments. 
So for the rest of the proof, we set 
$\vt=${\it ``left"}.  

$i)$ Assume otherwise, i.e., $(A)_\vt = KA$. 
Then for any $X\in M_n(K)$, we have $XA=rA$ for 
some $r\in K$. Consequently, each column of $A$ 
is a common eigenvector of all matrices $X\in M_n(K)$, 
which is clearly impossible unless the column 
is equal to zero. Therefore, we have $A=0$, 
which is a contradiction. 

$ii)$ Note first that since quasi-idempotents are preserved 
by taking conjugations, we may replace $A$ by any 
conjugation of $A$. 

Write $A=X\begin{pmatrix} I_k & 0 \\0  & 0 \end{pmatrix} Y$ 
for some $1\le k\le n$ and invertible $X, Y\in M_n(K)$. 
Replacing $A$ by $YAY^{-1}$, we have 
$A=YX\begin{pmatrix} I_k & 0 \\0  & 0 \end{pmatrix}$. 
Since $YX$ is invertible, the left ideal $(A)_\vt$ generated by $A$ is the same as the left ideal generated by 
$\begin{pmatrix} I_k & 0 \\0 & 0 \end{pmatrix}$.  
Hence, we may assume $A=\begin{pmatrix} I_k & 0 \\0 & 0 \end{pmatrix}$ for some $1\le k\le n$. 

Now, let $B=(b_{ij})\in M_n(K)$ such that $b_{i j}=1$ 
if $i=2$ and $j=1$; and $0$ otherwise. 
Then we have $B=BA\in (A)_\vt$. Since $B$ is nonzero and nilpotent, it follows from Lemma \ref{CyclicLemma} $i)$ 
that $B$ cannot be a quasi-idempotent, 
whence the statement follows. 
\epfv

\underline{\it Proof of Proposition \ref{ClsMin}:}\, 
The $(\Leftarrow)$ part follows directly from 
Proposition \ref{1D-MS}. To show the $(\Rightarrow)$ part,
we first show $\dim_K V=1$.

Assume otherwise. 
Then for any nonzero $A\in V$, the line $KA\ne V$ and hence, cannot be a \vms of $M_n(K)$, for $V$ is minimal. 
Applying Proposition \ref{1D-MS} to $A$, 
we see that $A$ must be a quasi-idempotent. 
Therefore, all elements of $V$ must be quasi-idempotents. 
Moreover, by Lemma \ref{CyclicLemma} $iii)$,  
for any nonzero $A\in V$, we have $(A)_\vt  \subset V$, 
whence all elements of $(A)_\vt$ are also quasi-idempotents. 
But this contradicts Lemma \ref{ClsMin-L} $ii)$. 

Now, write $V=KA$ for some $0\neq A\in M_n(K)$. 
Then from Proposition \ref{1D-MS} and 
Lemma \ref{ClsMin-L} $i)$, it follows that $A$ cannot 
be a quasi-idempotent.  
\epfv

Finally, we conclude this section with the following 
remarks on the maximal non-trivial  
$\vt$-Mathieu subspaces of $M_n(K)$.

In contrast to the minimal non-trivial  case,   
the situation for the maximal non-trivial  
\vmss of $M_n(K)$ becomes much more complicated. 
Even though Theorem \ref{Co-D1} classifies 
the co-dimension one maximal $\vt$-Mathieu subspace of $M_n(K)$, 
there are also many others (with different co-dimensions). 

For example, pick up any $A\in M_n(K)\backslash H$ 
(i.e., ${\rm Tr\,} A\ne 0$) such that $A$ is not a 
quasi-idempotent. Then by Proposition \ref{1D-MS}, 
the line $KA$ is a \vms of $\cA$. 
If $n\ge 2$, then by Proposition \ref{MaxMSs-Thm} 
or by counting dimensions, $KA$ is contained in 
at least one maximal non-trivial  \vms $W$ of $M_n(K)$. 
But, since $A\in W$ and $A\not \in H$, 
we have $W\ne H$. 

The situation for the two-sided case can be slightly improved by the following proposition.

\begin{propo}\label{FinalRmk} 
Let $V$ be a proper $K$-subspace of $M_n(K)$. 
Then $V$ is a Mathieu subspace of $M_n(K)$ 
iff $V$ does not contain any nonzero idempotent of $M_n(K)$. 
\end{propo}
\pf The $(\Leftarrow)$ part follows immediately from 
Lemma \ref{3Cases} and Corollary \ref{CharByIdem-cor1} 
since $\dim_K M_n(K)<\infty$. 

To show  $(\Rightarrow)$ part, assume otherwise, 
i.e., there exists a nonzero idempotent $A\in V$. 
Then by Lemma \ref{CyclicLemma} $iii)$, the ideal 
$(A)$ of $M_n(K)$ generated by $A$ 
is also contained in $V$. 
But, on the other hand, it is well-known that 
$M_n(K)$ is a simple $K$-algebra 
(e.g., see the lemma on p.\,$9$ in \cite{P}).
Hence $(A)=M_n(K)$, whence $V=M_n(K)$. 
But this contradicts the assumption that 
$V$ is proper. 
\epfv 

%assume that $V$ contains a nonzero idempotent $e \in M_n(K)$. Then by Lemma \ref{CyclicLemma} $iii)$, the (two-sided) ideal $(e)$ of $M_n(K)$ generated by $e$ is contained in $V$. But, one the other hand, it is well-known that $M_n(K)$ is a simple $K$-algebra  
%(e.g., see the lemma on p.\,$9$ in \cite{P}). So  
%we also have $(e)=M_n(K)$, whence $V=M_n(K)$.  
%But this contradict to the assumption 
%that $V$ is proper. \epfv 

It will be interesting if one can get 
a more explicit classification (other than the ones  
given by Theorem \ref{CharByIdem} and 
Proposition \ref{FinalRmk} plus the maximality) 
of all maximal non-trivial    
\vmss for matrix algebras $M_n(K)$ $(n\ge 2)$,   
or even more generally, for all finite dimensional or 
algebraic $K$-algebras.

\renewcommand{\theequation}{\thesection.\arabic{equation}}
\renewcommand{\therema}{\thesection.\arabic{rema}}
\setcounter{equation}{0}
\setcounter{rema}{0}

\section{\bf Strongly Simple Algebras} 
\label{S6}

As we have mentioned earlier, the notion of Mathieu subspaces can be viewed as a natural generalization of the notion of ideals.
Note that one of the most important families of 
(associative) algebras are simple algebras, i.e., 
the algebras that have no non-trivial  ideals. 
Then parallel to simple algebras, we have the 
following family of special algebras.

\begin{defi}\label{StrSimAlg}
Let $R$ be a commutative ring and $\cA$ a  $R$-algebra. 
We say that $\cA$ is a {\it strongly simple $R$-algebra}  
if $\cA$ has no non-trivial  $($two-sided$)$ Mathieu subspaces.
\end{defi}

Formally, one may also consider 
{\it left} (resp., {\it right}, {\it pre-two-sided}) 
{\it strongly simple algebras}, i.e., the algebras that have no non-trivial  
{\it left} (resp., {\it right}, {\it pre-two-sided}) Mathieu subspaces. 
But, as we will show in Theorem \ref{No-MS-Thm1} below, every (two-sided) 
strongly simple algebras is commutative. From this fact, it is easy to see that the notion of {\it left}, {\it right} or {\it pre-two-sided}) 
strongly simple algebras is actually equivalent to the notion of 
(two-sided) strongly simple algebras. In other words, an algebra 
is {\it left}, {\it right} or {\it pre-two-sided}  strongly simple 
iff it is (two-sided) strongly simple.

In this section, we give a characterization for strongly simple algebras $\cA$  over arbitrary commutative rings $R$. For convenience, 
throughout the rest of this section except 
in Corollary \ref{RingCase}, we assume that 
the base ring $R$ is contained in the 
$R$-algebra $\cA$. 
Note that by replacing $R$ 
by $R\cdot 1_\cA\subseteq \cA$, this condition will 
be satisfied. Furthermore, when $R$ is an integral domain, 
we denote by $K_R$ the field of 
fractions of $R$. Note that by 
{\it ``integral domains"} we always mean 
{\it commutative} domains. 

Under the assumption and notation above, the first main result of this 
section can be stated as follows.

\begin{theo}\label{No-MS-Thm1}
Let $R$ be a commutative ring and $\cA$ a  $R$-algebra. 
Then $\cA$ is a strongly simple $R$-algebra $($if and$)$ only if the following three statements hold: 
\begin{enumerate}
  \item[$i)$] $R$ is an integral domain; 
  \item[$ii)$] $\cA \simeq K_R$ as $R$-algebras;
  \item[$iii)$] $K_R$ as a  $R$-algebra is strongly simple.
  \end{enumerate}  
\end{theo}

One immediate consequence of the theorem above 
is the following corollary.

\begin{corol}\label{No-MS-Corol-1}
Let $R$ be a commutative ring and $\cA$ a  $R$-algebra. 
Assume that either $R$ is not an integral domain, or $\cA$ 
is not commutative, or $\cA$ 
is commutative but not a field. 
Then $\cA$ has at least one 
non-trivial Mathieu subspace. 
\end{corol}

%Note that the theorem also holds when the 
%``Mathieu subspace" is replace by ``$\vartheta$-Mathieu subspace" since any 
%non-trivial  Mathieu subspace of $\cA$ is a one-sided Mathieu subspace of $\cA$ 
%and the notions of one-sided Mathieu subspace and 
%(two-sided) Mathieu subspace coincide for the commutative 
%$R$-algebras $\cA=K_R$.   

In order to prove Theorem \ref{No-MS-Thm1}, 
we first need to show the following lemma 
which is the special case of the theorem 
when the base ring $R$ is a field $K$.

\begin{lemma}\label{FieldCase}
Let $K$ be a field and 
$\cA$ a  $K$-algebra. Then $\cA$ is 
strongly simple $($if and$)$ only if $\cA=K$.
\end{lemma}

\pf Assume otherwise, i.e., $\cA\ne K$. Then there exists $a\in \cA$ 
such that $a$ is linearly independent 
with $1\in \cA$ over $K$. Throughout the rest of the proof, 
we fix such an element $a$ and derive a contradiction 
as follows.

First, since every non-trivial  ideal of $\cA$ 
is a non-trivial  Mathieu subspace of $\cA$, we see that 
$\cA$ cannot have any non-trivial  ideals, which means 
that $\cA$ must be a simple $K$-algebra. 

Second, for any nonzero $b\in \cA$, 
the one-dimensional $K$-subspace 
$Kb\subset \cA$ is non-trivial but cannot be a 
Mathieu subspace of $\cA$. Then by Proposition \ref{1D-MS}, 
$b$ must be a quasi-idempotent of $\cA$. Therefore, 
all nonzero elements of $\cA$ are quasi-idempotents.

In particular, the element $a\in \cA$ fixed at the 
beginning is a quasi-idempotent. Replacing $a$ by   
a scalar multiple of $a$, we further  
assume from now on that $a$ is an idempotent 
which is linearly independent with $1\in\cA$. 

Next, with the two observations above in mind we consider the 
following two different cases. 

\underline{\it Case 1}\,: Assume $K\simeq \bZ_2$. Then in this case 
all elements of $\cA$ are actually idempotents instead of just being quasi-idempotents (since the only nonzero element of the base field $K$ 
is $1\in K$). It is well-known or from the simple argument below 
that $\cA$ in this case is actually a commutative algebra. 
Since $\cA$ is also simple, we see that $\cA$ in this case 
is actually a field extension of $\bZ_2$.   
 
Let $b, c\in \cA$. Then $b$, $c$ and $b+c$ 
are all idempotents. From the equations $(b+c)^2=b+c$; 
$b^2=b$ and $c^2=c$, it is easy to see that 
$bc=-cb=cb$. 

%Since $\cA$ is simple and commutative, we see that 
%all nonzero elements of $\cA$ must be invertible.

Now, let $a\in \cA$ be the idempotent fixed above.  
Since $a\ne 0$ and $\cA$ is a field, $a$ is invertible. 
Then by Lemma \ref{CyclicLemma} $i)$, we have 
$a\in K^\times$. But this contradicts 
our assumption that $a$ and $1$ are linearly 
independent over $K$. 

\underline{\it Case 2}\,: Assume $K\not \simeq \bZ_2$. 
Then there exists $r\in K^\times$ such that $r \ne -1$. 
Set $b\!:=1+r a$, where $a$ is as fixed before. 
Note that $b\ne 0$ since $1$ and $a$ are linearly 
independent over $K$. Then we have $b^2=s b$ 
for some $s \in K^\times$. More precisely, 
we have 
\begin{align*}
s(1+r a)&=(1+ r a)^2=
1+2r a+r^2 a^2 \\
&=1+2r a +r^2 a=
1+(2+r)r a.
\end{align*}
By comparing the coefficients of 
$1$ and $a$ in the equation above, 
we get
\begin{align*}
\begin{cases}
s =1,\\
sr=(2+r)r .
\end{cases}
\end{align*}
Solving the equation above, we get $r=-1$, 
which is a contradiction again. 
Therefore, the lemma holds.
\epfv

The following lemma will also be important to us.

\begin{lemma}\label{GoingUp}
Let $S$ be a subring of a  $R$-algebra $\cA$ such that 
$R\subseteq S\subseteq Z(\cA)$, where $Z(\cA)$ 
denotes the center of $\cA$.
Assume that $\cA$ as a  $R$-algebra is strongly simple. 
Then $\cA$ as a  $S$-algebra is also strongly simple.
\end{lemma}
\pf Since $S\subseteq Z(\cA)$, 
$\cA$ can also be viewed as a  $S$-algebra (in 
the obvious way). Moreover, since $R\subseteq S$, 
every $S$-subspace of $\cA$ is also a $R$-subspace of $\cA$. 
With these observations, the lemma follows 
immediately from the definition of Mathieu subspaces 
(see Definitions \ref{Def-MS} and \ref{Def-MS4}) 
and that of strongly simple algebras 
(see Definition \ref{StrSimAlg}).
\epfv

Now we can prove Theorem \ref{No-MS-Thm1} 
as follows.\\

\underline{\it Proof of Theorem \ref{No-MS-Thm1}:} \, 
First, let $0\ne r\in R\subseteq \cA$. Since $r$ commutes with all elements of $\cA$, $\cA r$ is a nonzero 
(two-sided) ideal and hence, also a nonzero 
Mathieu subspace of $\cA$. 
Since $\cA$ is a strongly simple $R$-algebra, 
we have $\cA r=\cA$. In particular, 
$1\in \cA r$ and $r$ is invertible 
in $\cA$. 
Therefore, all nonzero elements of $R$ 
are invertible in $\cA$, whence $R$ must 
be an integral domain, i.e., 
the statement $i)$ in 
the theorem holds.

Furthermore, since $R\subseteq \cA$, 
we may also assume that $\cA$ 
contains the field of fractions $K_R$ of $R$. 
Since all elements of $R$ are central elements 
of $\cA$, it is easy to check that so are 
all elements of $K_R\subseteq \cA$. Therefore,  
$\cA$ can also be viewed as a $K_R$-algebra. 

Now, since $\cA$ is strongly simple 
as a  $R$-algebra, by 
Lemma \ref{GoingUp} with $S=K_R$, 
it is also strongly simple 
as a  $K_R$-algebra. Then 
by Lemma \ref{FieldCase}, 
we have $\cA=K_R$. Therefore, 
the statement $ii)$ in the theorem 
holds. The statement $iii)$ follows 
from the statement $ii)$ and 
our assumption on the  
$R$-algebra $\cA$.  
\epfv

From Theorem \ref{No-MS-Thm1}, we see that 
in order to classify all strongly simple algebras, 
it suffices to classify all the integral domains 
$R$ whose field of fractions $K_R\ne R$ and 
as a $R$-algebra is strongly simple. 

We have not succeeded in classifying 
this special family of integral domains. 
Instead, we show next that no 
Noetherian domain or Krull domain belongs 
to this family. To do so, we first need 
to prove the following lemma. 
   
\begin{lemma}\label{ValuationLemma}
Let $R$ be an integral domain with $R\ne K_R$. 
Assume that there exists a non-trivial real-valued 
additive valuation $\nu$ of $K_R$ such that 
$\nu(r)\ge 0$ for all $r\in R$. 
Then $K_R$ as a  $R$-algebra 
is not strongly simple.
\end{lemma}
\pf Since $\nu$ is non-trivial, i.e., 
$\nu(a)\ne 0$ for some $0\ne a\in K_R$, there exists 
a positive $\beta\in \bR$ such that 
$M_\beta\!:=\{a\in K_R\,|\, \nu(a)\ge \beta\}\ne 
0$. Note that $M_\beta\ne K_R$ either 
since for each $a\in M_\beta$, we have 
$\nu(a^{-1})=-\nu(a)<0$, 
whence $a^{-1}\not \in M_\beta$. 
 
Furthermore, by our assumption that $\nu(r)\ge 0$ 
for all $r\in R$, it is easy to check that $M_\beta$ is a 
Mathieu subspace of $K_R$. Therefore, $M_\beta$ is a 
non-trivial  Mathieu subspace of $K_R$, whence $K_R$ 
is not a strongly simple $R$-algebra.
\epfv

For general discussions on valuations, and also on 
Krull domains needed below, see \cite{Sc}, 
\cite{R}, \cite{AM}, \cite{ZS}, \cite{Bour} and \cite{Fo}.  

\begin{propo}\label{SpecialDomains}
Let $R$ be a Krull domain or a Noetherian domain 
such that $R\ne K_R$, i.e., $R$ is not a field. 
Then no $R$-algebra is strongly simple. 
Equivalently, every $R$-algebra $\cA$ 
has at least one non-trivial  Mathieu subspace. 
\end{propo}

\pf Note that by Theorem \ref{No-MS-Thm1}, 
it suffices to show that $K_R$ as a  $R$-algebra 
is not strongly simple. 

Assume first that $R$ is a Krull domain.  
Since $R$ is not a field, by the very definition 
of Krull domains (e.g., see p.\,$480$ in \cite{Bour}), 
we see that $R$ satisfies the 
hypothesis in Lemma \ref{ValuationLemma}. Hence,  
by Lemma \ref{ValuationLemma} $K_R$ cannot be 
a strongly simple $R$-algebra. 

Now, assume that $R$ is a Noetherian domain. 
Let $\bar R$ be the integral closure of $R$ 
in $K_R$. Then by the Mori-Nagata integral 
closure theorem (see Theorem 4.3, p.\,18 in \cite{Fo} or  
Corollary 2.3, p.\,161 in \cite{Hu}), 
$\bar R$ is a Krull domain. Note that since $R$ is not 
a field, it is well-known (e.g., see Proposition $5.7$, 
p.\,$61$ in \cite{AM}) that $\bar R$ is not a field either.      

Furthermore, since the field of fractions 
$K_{\bar R}$ of $\bar R$ is the same as $K_R$, 
by the Krull domain case that we just 
proved above, $K_R$ is not strongly simple 
as a  $\bar R$-algebra, and by Lemma \ref{GoingUp} 
with $S=\bar R$, $K_R$ is not 
strongly simple as a  $R$-algebra either. 
\epfv

Since $\bZ$ and all its quotient rings are obviously Noetherian, 
from Theorem \ref{No-MS-Thm1} and Proposition 
\ref{SpecialDomains} we immediately have the following 
classification for {\it strongly simple rings} $\cA$, 
i.e., strongly simple algebras $\cA$ over $\bZ$ 
(without the convenient assumption $\bZ\subseteq \cA$). 

\begin{corol}\label{RingCase} 
Let $\cA$ be an arbitrary  
commutative or noncommutative ring.  
Then $\cA$ as a $\bZ$-algebra is strongly simple 
iff $\cA\simeq \bZ_p$ for some prime $p>0$.  
In other words, all rings $($as $\bZ$-algebras$)$  
except the finite fields $\bZ_p$'s have 
non-trivial Mathieu subspaces.   
\end{corol}

Next, we conclude this section with 
the following remarks. 

\begin{rmk}
$i)$ By Lemma \ref{GoingUp}, we see that 
Proposition \ref{SpecialDomains} also holds 
if there exists a Noetherian or Krull domain $S$ 
of $K_R$ such that $S$ is not a field and 
$S$ contains $R$.    

$ii)$ 
After an earlier version of this paper was circulated,  
M. de Bondt \cite{Bon} has  
recently found some examples of integral domains $R$ such that $R$ 
is not a field and $K_R$ is strongly simple as a  $R$-algebra. 
He also showed that for any integral domain $R$ that has at least 
one prime ideal of height one, the field of fractions $K_R$ 
as a $R$-algebra is not strongly simple. 
Therefore, by Theorem \ref{No-MS-Thm1} we see that 
Proposition \ref{SpecialDomains} actually holds for  
all the integral domains with prime ideals 
of height one. 
\end{rmk}

%It seems plausible that Corollary \ref{SpecialDomains} 
%should hold for all integral domains. So we formulate 
%it here as a conjecture.  
%
%\begin{conj}
%For every integral domain $R$ that is not a field, its field 
%of fractions $K_R$ as a  $R$-algebra is not strongly simple.
%\end{conj}

%Evidence: Exercise 5 on p.\,450 in \cite{Bour}.

\renewcommand{\theequation}{\thesection.\arabic{equation}}
\renewcommand{\therema}{\thesection.\arabic{rema}}
\setcounter{equation}{0}
\setcounter{rema}{0}

\section{\bf Quasi-Stable Algebras}\label{S7}
First, let's introduce the following notions    
for associative algebras. 

\begin{defi}\label{q-StaAlg}
Let $\cA$ be an associative $R$-algebra. 
We say that $\cA$ is 
{\it $\vartheta$-quasi-stable} 
$($resp., $\vartheta$-stable$)$ if every $R$-subspace 
$V$ of $\cA$ with $1\not \in V$ is a 
\vms $($resp., $\vartheta$-ideal$)$ of $\cA$.
\end{defi}

For the justifications of the terminologies 
in the definition above, see Section $3$ 
in \cite{GMS}. 

In contrast to {\it strongly simple algebras} studied 
in the previous section, which 
have as less \vmss as possible, 
{\it $\vartheta$-quasi-stable algebras} 
by Corollary \ref{OneLemma} 
are the algebras that have as many \vmss as possible. 
One of the motivations for the study of 
$\vartheta$-quasi-stable algebras comes 
from the following proposition and the 
corollary followed.

\begin{propo}\label{MotivStable}
Let $\cA$ and $\cB$ be $R$-algebras and $\phi:\cB\to \cA$ 
a  $R$-algebra homomorphism. Assume that $\cA$ is 
$\vartheta$-quasi-stable. Then for every $R$-subspace 
$V$ of $\cA$ such that $1_\cA \not \in V$, 
the pre-image $\phi^{-1}(V)$ is a   
\vms of $\cB$.
\end{propo}
\pf Since $1_{\cA}\not \in V$ and $\cA$ is a 
$\vartheta$-quasi-stable $R$-algebra, 
we have that $V$ is a \vms of $\cA$. 
Then by Proposition \ref{Pull-Back},  
$\phi^{-1}(V)$ is a \vms of $\cB$.  
\epfv

\begin{corol}\label{MotivStable-Corol}
Let $\cB$ be a  $R$-algebra and $I$ an ideal 
of $\cB$ such that $\cB/I$ 
is a $\vartheta$-quasi-stable 
$R$-algebra. Then every $R$-subspace $M$ of \, $\cB$ 
with $I\subseteq M$ and $1\not \in M$ 
is a \vms of $\cB$.
\end{corol}
\pf If $I=\cB$, the corollary holds  
vacuously. So we assume $I\ne \cB$. 
Let $\cA\!:=\cB/I$ 
and $\pi: \cB\to \cA$ the quotient $R$-algebra 
homomorphism. Set $V\!:=\pi(M)$. 
Then by the assumptions $1_\cB\not \in M$  
and $I\subseteq M$, it is easy to check that 
$1_{\cA}\not \in V$ and $M=\pi^{-1}(V)$. 
Applying Proposition \ref{MotivStable} to 
the $R$-subspace $V\subset \cA$ with 
$\phi=\pi$, we see that the corollary follows.  
\epfv

One family of quasi-stable 
$R$-algebras is given by the following 
proposition.

\begin{propo}\label{q-Quasi/R}
Let $\cA$ be a $R$-algebra such that $\cA$ is integral 
over $R$ and every element of $\cA$ is either invertible or 
nilpotent. Then $\cA$ is a $\vt$-quasi-stable $R$-algebra.  
\end{propo}
\pf Let $V$ be a  $R$-subspace of $\cA$ such that 
$1\not \in V$. Since $\cA$ is integral over $R$, 
by Lemma \ref{Inv-Alg-L1-New} the radical $\sqrt V$ 
of $V$ does not contain any invertible element of $\cA$. 
Hence by our assumption on $\cA$, we have  
$\sqrt V\subseteq {\rm nil\,}(\cA)$. 
Then by Lemma \ref{radical-Lemma}, 
$V$ is a \vms of $\cA$. 
Hence the proposition follows. 
\epfv

\begin{corol}\label{ArtinLocal}
Every left or right Artinian local $R$-algebra 
$\cA$ that is integral over $R$ is 
$\vt$-quasi-stable. 
In particular, every commutative 
Artinian local ring as a $\bZ$-algebra is quasi-stable  
if it is integral over $\bZ$.  
\end{corol} 
\pf Since $\cA$ is local, it's Jacobson radical 
${\bf J}(\cA)$ is also the unique maximal left 
ideal of $\cA$. Hence, all non-invertible elements 
of $\cA$ are contained in ${\bf J}(\cA)$. 
Since $\cA$ is left or right Artinian, it is well-known 
(e.g., see the proposition on p.\,$61$ in \cite{P})
that the Jacobson radical ${\bf J}(\cA)$ is nilpotent, 
i.e., ${\bf J}(\cA)^k=0$ for some $k\ge 1$. 
Consequently, all the elements in ${\bf J}(\cA)$ 
are nilpotent. Therefore, all elements of $\cA$ are either 
invertible or nilpotent, and by Proposition \ref{q-Quasi/R}, 
$\cA$ is $\vt$-quasi-stable.  
\epfv

Next, we give the following classification 
for $\vartheta$-quasi-stable algebras $\cA$ 
over arbitrary fields $K$. 

\begin{theo}\label{Class-Qstable}
Let $K$ be a field and $\cA$ a  $K$-algebra.  
Then $\cA$ is $\vartheta$-quasi-stable iff either 
$\cA \simeq K \dot{+} K$ or $\cA$ is 
an algebraic local $K$-algebra. 
%that satisfies one of the following three equivalent 
%statements:
% 
%$1)$ $\cA$ is a local $K$-algebra;
%
%$2)$ $\cA$ has no non-trivial idempotent;
%
%$3)$ every element of $\cA$ is 
%either nilpotent or invertible. 
\end{theo}

Two remarks on the theorem above are as follows. 

First, by Corollary \ref{Alg-3Equiv} we see that for any algebraic 
$K$-algebra $\cA$, $\cA$ is local 
iff every element of $\cA$ is either nilpotent or invertible. Therefore, by Theorem \ref{Class-Qstable} 
we see that Proposition \ref{q-Quasi/R} with $R=K$ 
actually has covered most of the $\vartheta$-quasi-stable 
algebras over $K$. 
 
Second, from Theorem \ref{Class-Qstable}, 
Corollary \ref{Alg-3Equiv}, Lemma \ref{3Cases}, 
Corollary \ref{CharByIdem-cor1}, 
or from the proof of Theorem \ref{Class-Qstable}  
given below, it follows that the $\vartheta$-quasi-stableness 
for algebras over fields actually does not depend on the specifications of 
$\vartheta$. More precisely, we have the following corollary.

\begin{corol}
Let $K$ be a field and $\cA$ a  $K$-algebra. 
Then $\cA$ is $\vt$-quasi-stable for one 
specification of $\vt$ iff $\cA$ is $\vt$-quasi-stable for all  
specifications of $\vt$ iff $\cA$ is $($two-sided$)$ quasi-stable.
\end{corol}

To prove Theorem \ref{Class-Qstable}, we start with the following lemma. 

\begin{lemma}\label{Stabe2Alg}
Every $\vartheta$-quasi-stable $K$-algebra $\cA$   
is algebraic over $K$.
\end{lemma}
\pf Assume otherwise and let $a$ be a (nonzero) element of 
$\cA$ which is transcendental over $K$. Denote by $V$ 
the $K$-subspace of $\cA$ spanned 
by $a^{2k}$ $(k\ge 1)$ over $K$. 
Then we have $1\not \in V$, otherwise $a$ would 
be algebraic over $K$. So $V$ is a \vms of $\cA$,   
for $\cA$ is $\vartheta$-quasi-stable. 
 
Since $(a^2)^m=a^{2m}\in V$ for all $m\ge 1$, there exists 
a large enough $N\ge 1$ such that $a^{2N+1}=(a^2)^N a \in V$.  
But this means that the odd power $a^{2N+1}$ 
can be written as a linear combination of 
some even powers of $a$, 
whence $a$ is algebraic over $K$. 
Hence, we get a contradiction. 
\epfv 

%
%\begin{lemma}\label{Nil-Inv-Case}
%Let $\cA$ be any algebraic $R$-algebra such that 
%the condition $(2)$, or equivalently, the condition 
%$(3)$ in Theorem \ref{Class-Qstable} holds. 
%Then $\cA$ is quasi-stable.  
%\end{lemma}  
%
%\pf Let $V$ be any $R$-subspace of $\cA$ such that 
%$1\not\in V$. Note that for any $a\in \sqrt V$, 
%by Lemma \ref{Inv-Alg-L1}, $a$ cannot be invertible in $\cA$, 
%whence by the condition $(3)$ in Theorem \ref{Class-Qstable}, 
%$a$ must be nilpotent. Therefore, we have $\sqrt V \subseteq \nil (\cA)$.  
%Then by Lemma \ref{radical-Lemma}, the proposition holds.  
%\epfv

\underline{\it Proof of Theorem \ref{Class-Qstable}:} \, 
$(\Leftarrow)$ Assume first that $\cA\simeq K\dot{+}K$, 
then it is easy to check that the only non-trivial  
idempotents of $\cA$ are $a\!:=(1, 0)$ and 
$b\!:=(0, 1)$. Note that the lines $Ka$ and $Kb$ 
are obviously ideals of $\cA$ and hence, also 
Mathieu subspaces of $\cA$. Then by Proposition \ref{1D-MS}, 
it is easy to see that every non-trivial  subspace $V$ of $\cA$ 
(which is necessarily a line of $\cA$) with 
$1_{\cA}=(1, 1)\not\in V$ is a Mathieu 
subspace of $\cA$. Therefore,  
$\cA$ is quasi-stable and hence,  
also $\vt$-quasi-stable 
for all possible $\vt$.  

Now assume that $\cA$ is an algebraic local $K$-algebra. 
Then by Corollary \ref{Alg-3Equiv}, 
$\cA$ has no non-trivial idempotent. 
Let $V$ be a  $R$-subspace of $\cA$ such that 
$1\not\in V$. Then $V$ contains no nonzero 
idempotent of $\cA$. By Lemma \ref{3Cases} and 
Corollary \ref{CharByIdem-cor1}, $V$ is a \vms of $\cA$. 
Therefore, $\cA$ is $\vartheta$-quasi-stable.

%Let $V$ be any $R$-subspace of $\cA$ such that 
%$1\not\in V$. Then for any $a\in \sqrt V$, 
%by Lemma \ref{Inv-Alg-L1}, $a\not \in \cA^\times$, 
%whence $a$ must be nilpotent. Therefore, 
%we have $\sqrt V \subseteq \nil (\cA)$. 
%Then by Lemma \ref{radical-Lemma}, 
%$V$ is a Mathieu subspace and hence, 
%also a $\vartheta$-Mathieu 
%subspace of $\cA$. Therefore, 
%$\cA$ is $\vartheta$-quasi-stable.

$(\Rightarrow)$ Assume that $\cA$ is not an algebraic local 
$K$-algebra. Then by Corollary \ref{Alg-3Equiv}, 
$\cA$ has at least one non-trivial idempotent, say, $e\in \cA$. 
Note that $e$ is linearly independent with $1\in \cA$ 
over $K$ since the only idempotents of $K$ 
are $0, 1\in K$. 

Let $\cB$ be the two-dimensional $K$-subspace 
of $\cA$ spanned by $1, e\in \cA$ over $K$. Then 
it is easy to check that $\cB$ is actually a  $K$-subalgebra 
of $\cA$ which is isomorphic to the $K$-algebra $K\dot{+}K$ 
via the following $K$-algebra isomorphism:  
\begin{align*}
\phi:K\dot{+}K\quad &\longrightarrow \quad \qquad \, \cB\\ 
(r,\,\, s) \quad &\longrightarrow \quad r(1-e)+ se.
\end{align*}

Next we show $\cB=\cA$, from which the theorem 
will follow. 

First, by the fact that $1$ and $e$ are linearly 
independent over $K$, we have $1\not \in K(1-e)$ and $1\not \in Ke$. 
Second, since $\cA$ is $\vartheta$-quasi-stable, both 
$Ke$ and $K(1-e)$ are \vmss of $\cA$. But, on the other hand, 
since $e$ and $(1-e)$ are non-trivial idempotents, 
it follows from Proposition \ref{1D-MS} 
that $Ke$ and $K(1-e)$ are actually 
$\vartheta$-ideals of $\cA$. 

Assume $\vartheta=${\it ``left"}, {\it ``pre-two-sided\hspace{.1mm}"} or 
{\it ``two-sided\hspace{.1mm}"}. Then for each $a\in \cA$, 
we have 
\begin{align}
ae&=r e,\\
a(1-e)&=s (1-e),
\end{align}
for some $r, s\in K$. 

Taking the sum of the two equations above, 
we get $a=r e+s(1-e)$, whence $a\in \cB$. 
Therefore, we do have $\cB=\cA$ 
when $\vartheta\ne${\it ``right"}. 
The case $\vartheta=${\it ``right"} 
can be proved similarly. Therefore, 
the theorem holds. 
\epfv

From Theorem \ref{Class-Qstable}, 
we immediately have the following 
examples of quasi-stable $K$-algebras.

\begin{exam}\label{Artin-Exam}
$1)$ every algebraic field extension of $K$ 
or more generally, every algebraic division algebra over 
$K$ is a quasi-stable $K$-algebra.

$2)$ Let $p$ be a prime and $\cA\!:=\bZ/(p^k)$ 
for some $k\ge 1$. Then $\cA$ as a $\bZ$-algebra 
is algebraic and local and hence, a quasi-stable 
$\bZ$-algebra. Actually, $\cA$ 
is also a stable $\bZ$-algebra since 
every $\bZ$-subspace of $\cA$ is an 
ideal of $\cA$.   

$3)$ Let $K$ be a field and $t$ a free variable. 
For every $k\ge 1$ and irreducible $f(t)\in K[t]$, 
the quotient algebra $\cA\!:=K[t]/(f^k)$ is an algebraic 
and local $K$-algebra. Therefore, 
$\cA$ is a quasi-stable $K$-algebra.   
\end{exam}

Note that all the quasi-stable algebras 
in the example above are Artinian. 
However, this is not always the case. 

\begin{exam}\label{No-NeoExam}
Let $\cB=K[x_i\,|\, i\ge 1]$ 
be the polynomial algebra over $K$ in the infinitely many commutative 
free variables $x_i$ $(i\ge 1)$, and $I$ 
the ideal of $\cB$ generated by $x_i^{i+1}$ $(i\ge 1)$.
Set $\cA\!:=\cB/I$. Then it is easy to see that $\cA$ 
is an algebraic local $K$-algebra whose  
maximal   ideal $\frak m$ is the ideal generated 
by the images of $x_i$ $(i\ge 1)$ in $\cA$. 
Hence, $\cA$ by Theorem \ref{Class-Qstable} is a 
quasi-stable $K$-algebra. 

On the other hand, since the maximal ideal $\frak m$ of 
$\cA$ is obviously not finitely generated, $\cA$ is not 
Noetherian and hence, not Artinian either. 
\end{exam}

The following proposition generalizes the construction in 
Example \ref{Artin-Exam}, $2)$ and $3)$ for 
quasi-stable algebras.

\begin{propo}\label{Noe2Qstable}
Let $\cA$ be a commutative 
$K$-algebra 
and $\frak{m}$ a maximal ideal of $\cA$ such that 
$\cA/\frak {m}$ is an algebraic field 
extension of $K$. Then for every $k\ge 1$, 
$\cA/\frak {m}^k$ is a quasi-stable 
$K$-algebra.  
\end{propo}

\pf It is easy to see that $\cA/\frak {m}^k$ 
is a local $K$-algebra with the maximal ideal 
$\frak m / \frak {m}^k$. Then by Theorem \ref{Class-Qstable}, 
we only need to show that $\cA/\frak {m}^k$ 
is algebraic over $K$. 

Let $a\in \cA$. Since $\cA/\frak {m}$ 
is algebraic over $K$, it is easy to see 
that there exists a nonzero polynomial 
$f(t)\in K[t]$ such that $f(a)\in \frak m$.   
Then we have $f^k(a)\in \frak m^k$ and 
$f^k(\bar a)=0$, where $\bar a$ denotes the image of 
$a$ in $\cA/\frak {m}^k$. Therefore $\bar a$ is 
algebraic over $K$ for all $a\in \cA$,  
whence $\cA/\frak {m}^k$ is algebraic over $K$.   
\epfv

From Proposition \ref{Noe2Qstable} and 
Corollary \ref{MotivStable-Corol}, 
we immediately have the following corollary.

\begin{corol}\label{App2Noe}
Let $\cA$ and $\frak{m}$ be as 
Proposition \ref{Noe2Qstable} and   
$V$ a  $K$-subspace of $\cA$. 
Assume that $1\not \in V$ and $\frak{m}^k\subseteq V$ 
for some $k\ge 1$. 
Then $V$ is a Mathieu subspace of $\cA$.
\end{corol}
%\pf The corollary follows immediately from 
%Lemma \ref{Noe2Qstable} and Proposition \ref{MotivStable} 
%with $\phi=\pi$, where $\pi: \cA\to \cA/\frak{m}^k$ is 
%the quotient homomorphism of $K$-algebras. 
%\epfv

In contrast to $\vartheta$-quasi-stable $K$-algebras, 
$\vartheta$-stable $K$-algebras do not seem very 
interesting. But, for the completeness and also 
for the purpose of comparison with 
$\vartheta$-quasi-stable algebras, here we conclude 
this paper with the following classification of 
$\vartheta$-stable $K$-algebras. 

\begin{propo}\label{ClassiStable}
Let $K$ be a field and $\cA$ a  $K$-algebra.  
Then $\cA$ is $\vartheta$-stable iff 
one of the following two statements holds:
\begin{enumerate}
  \item[$1)$] $\cA=K$;
  \item[$2)$] $K\simeq \bZ_2$ and $\cA\simeq \bZ_2\dot{+}\bZ_2$.
  \end{enumerate}    
\end{propo}
\pf The $(\Leftarrow)$ part of the proposition can be easily checked. 
To show the $(\Rightarrow)$ part, we assume $\cA\ne K$, 
and claim first that the following equation holds: 
\begin{align}\label{ClassiStable-pe1}
\cA^\times =K^\times.
\end{align} 

Assume otherwise and let 
$a\in \cA^\times \backslash K$. 
Then $1\not \in Ka$.   
Since $\cA$ is $\vartheta$-stable, 
$Ka$ is a $\vartheta$-ideal of $\cA$. 
But for any $\vt$, this implies $1=a^{-1}a=aa^{-1}\in Ka$, 
which is a contradiction. 
Therefore, Eq.\,(\ref{ClassiStable-pe1}) 
does hold.

On the other hand, since every $\vartheta$-stable algebra 
is obviously $\vartheta$-quasi-stable, hence $\cA$ by our 
hypothesis is also $\vartheta$-quasi-stable. 
Then by Theorem \ref{Class-Qstable}, we have that either  
$\cA\simeq K\dot{+}K$ or $\cA$ is an algebraic 
local $K$-algebra.

In the latter case, it follows from Corollary \ref{Alg-3Equiv} 
that all elements of $\cA$ are either nilpotent or invertible. 
Then by Eq.\,(\ref{ClassiStable-pe1}), all elements in 
$\cA\backslash K$ are nilpotent. But this is impossible 
by the argument below.

Let $a\in \cA\backslash K$ and set $b\!:=1-a$. 
Then $b\not \in K$. Hence, both $a$ and $b$ 
are nilpotent. But, on the other hand, since $a$ 
is nilpotent, $b$ has inverse $\sum_{i\ge 0} a^i$ 
in $\cA$. Therefore, we have $b\in \cA^\times$ (and  
$b\not \in K$), which contradicts  
Eq.\,(\ref{ClassiStable-pe1}).  

Therefore, we must have $\cA\simeq K\dot{+}K$.
If $K\not \simeq \bZ_2$, then there exist $r, s\in K^\times$ 
such that $r\ne s$. Set $a\!:=(r, s)$. Then $a\in \cA^\times$ 
and $a$ does not lie in the base field 
$K\simeq K\cdot 1_\cA\subset \cA$, 
since $1_\cA=(1, 1)$.  But this contradicts 
Eq.\,(\ref{ClassiStable-pe1}) again. 
Hence, the theorem follows. 
\epfv

{\bf Acknowledgments}\,\,\,\,The author is very grateful to Michiel de Bondt for pointing out some mistakes and misprints in the earlier version of the paper, for making many valuable suggestions,   
and also for sending the author his recent results 
on strongly simple algebras discussed in Section \ref{S6}.

%{\small \sc Department of Mathematics, Illinois State University,
%Normal, IL 61790-4520.}
%
%{\em E-mail}: wzhao@ilstu.edu.
%
%\end{document}

\end{document}